\documentclass{amsart}
\usepackage{latexsym}
\def\qed{\hfill $\Box$}
\usepackage{amsfonts,amssymb,amsmath}
\usepackage{amsthm}
\usepackage{graphicx,color}
\usepackage{color}
\usepackage{amsmath} 
\usepackage{multirow,bigdelim}
\newtheorem{theorem}{Theorem}	
\newtheorem{lemma}{Lemma}
\newtheorem{corollary}{Corollary}		
\newtheorem{proposition}{Proposition}

\newtheorem{definition}{Definition}

\newtheorem{remark}{Remark}

\newfont{\bg}{cmr9 scaled\magstep2}
\newcommand{\bigzerol}{\smash{\lower1.0ex\hbox{\bg 0}}}

\newcommand{\R}{\mathbb{R}}
\def\<{\color{black}} 
\def\>{\color{black}} 

\makeatletter
\newcommand{\tpitchfork}{%
  \vbox{
    \baselineskip\z@skip
    \lineskip-.52ex
    \lineskiplimit\maxdimen
    \m@th
    \ialign{##\crcr\hidewidth\smash{$-$}\hidewidth\crcr$\pitchfork$\crcr}
  }%
}
\makeatother
\begin{document}

\title[Composing linearly perturbed mappings 
and immersions/injections]{Composing generic linearly perturbed mappings 
and immersions/injections}


\author{Shunsuke Ichiki}

\dedicatory{{\it Dedicated to Professor Takashi Nishimura on the occasion of 
his 60th birthday}}
\address{
Graduate School of Environment and Information Sciences,  
Yokohama National University, 
Yokohama 240-8501, Japan}
\email{ichiki-shunsuke-jb@ynu.jp}





\begin{abstract}
Let $N$ (resp., $U$) be a manifold (resp., an open subset of $\R^m$). 
Let $f:N\to U$ and $F:U\to \R^\ell$ be an immersion and a $C^{\infty}$ mapping, respectively.  
Generally, the composition $F\circ f$ does not necessarily 
yield a mapping transverse  
to a given subfiber-bundle of $J^1(N,\mathbb{R}^\ell)$. 
Nevertheless, in this paper, for any $\mathcal{A}^1$-invariant fiber, we show that 
composing generic linearly perturbed mappings of $F$ and the given immersion $f$ yields a mapping transverse to 
the subfiber-bundle of $J^1(N,\mathbb{R}^\ell)$ with the given fiber. 
Moreover, we show a specialized transversality theorem on crossings 
of compositions of 
generic linearly perturbed mappings of a given mapping $F:U\to \mathbb{R}^\ell$ 
and a given injection $f:N\to U$. 
Furthermore, applications of the two main theorems are given.
\end{abstract}

\maketitle

\section{Introduction}\label{section 1}
Throughout this paper, let $\ell$, $m$ and $n$ stand for positive integers. 
In this paper, unless otherwise stated, all manifolds and mappings belong to class $C^{\infty}$ 
and all manifolds are without boundary. 
Let $\pi:\mathbb{R}^m\to \mathbb{R}^\ell$, $U$ and 
$F:U\to \mathbb{R}^\ell$ be a linear mapping, 
an open subset of $\R^m$ 
and a mapping, respectively. 

Set 
\begin{eqnarray*} 
F_\pi=F+\pi. 
\end{eqnarray*} 
Here, the mapping $\pi$ in $F_\pi=F+\pi$ is restricted to $U$. 

Let $\mathcal{L}(\mathbb{R}^{m},\mathbb{R}^{\ell})$ 
be the space consisting of all linear mappings 
of $\mathbb{R}^{m}$ into $\mathbb{R}^{\ell}$. 
Remark that \<we have \> the natural identification 
$\mathcal{L}(\mathbb{R}^{m},\mathbb{R}^{\ell})=(\mathbb{R}^m)^\ell$.  
An $n$-dimensional manifold is denoted by $N$. 
For 
a given mapping $f:N \to U$, 
a property of mappings $F_\pi \circ f:N\to \mathbb{R}^\ell$ will be said to be 
true for a {\it generic mapping} 
if there exists a subset $\Sigma$ with Lebesgue measure zero of 
$\mathcal{L}(\mathbb{R}^{m},\mathbb{R}^{\ell})$ 
such that for any $\pi \in \mathcal{L}(\mathbb{R}^{m},\mathbb{R}^{\ell})-\Sigma$, 
the mapping $F_\pi \circ f:N\to \mathbb{R}^{\ell}$ has the property. 
In the case $F=0$, 
by John Mather, for a given embedding $f:N\to \mathbb{R}^m$, 
a generic mapping $\pi\circ f:N\to \mathbb{R}^\ell$ $(m>\ell)$ 
is investigated in the celebrated paper \cite{GP}. 
The main theorem in \cite{GP} yields many applications. 
On the other hand, 
in this paper, for a given immersion or a given injection 
$f:N\to U$, 
 a generic mapping   
$F_\pi \circ f:N\to \mathbb{R}^\ell$ 
is investigated, where $\ell$ is an arbitrary positive integer which 
may possibly satisfy $m\leq \ell$. 

The main purpose of this paper is to show two main theorems 
(Theorems \ref{main} and \ref{main2} in Section \ref{section 2}) 
and to give some of their applications. 
The first main theorem 
(Theorem \ref{main}) is as follows. 
Let $f:N\to U$ (resp., $F:U\to \R^\ell$) be an immersion 
(resp., a mapping). 
Then, 
generally, the composition $F\circ f$ does not necessarily 
yield a mapping transverse 
to a given subfiber-bundle of the jet bundle $J^1(N,\mathbb{R}^\ell)$. 
Nevertheless, Theorem \ref{main} asserts that  
for any $\mathcal{A}^1$-invariant fiber, a 
generic mapping $F_\pi \circ f$ yields a mapping transverse to 
the subfiber-bundle of $J^1(N,\mathbb{R}^\ell)$ with the given fiber. 
The second main theorem (Theorem \ref{main2}) is a specialized transversality theorem on crossings of 
a generic mapping $F_\pi \circ f$, where $f:N\to U$ is a given injection and 
$F:U\to \mathbb{R}^\ell$ is a given mapping. 


For a given immersion (resp., injection) $f:N\to U$, 
the following (1)-(4) (resp., (5)) are obtained as applications of Theorem \ref{main} (resp., Theorem \ref{main2}). 

\begin{enumerate}
\item[(1)] If $(n,\ell)=(n,1)$, 
then a generic function $F_{\pi}\circ f:N\to \mathbb{R}$ is a Morse function. 
\item[(2)] If $(n,\ell)=(n,2n-1)$ and $n\geq 2$, 
then any singular point of a generic mapping $F_{\pi}\circ f:N\to \mathbb{R}^{2n-1}$ 
is a singular point of Whitney umbrella. 
\item[(3)]If $\ell \geq 2n$, 
then a generic mapping $F_{\pi}\circ f:N\to \mathbb{R}^{\ell}$ is an immersion.  
\item[(4)] A generic mapping $F_{\pi}\circ f:N\to \mathbb{R}^{\ell}$ 
has corank at most $k$ singular points 
\<(for the definition of corank at most $k$ singular points, see Subsection \ref{application1})\>, 
where $k$ is the maximum integer satisfying 
$(n-v+k)(\ell-v+k)\leq n$ $(v={\rm min}\{n,\ell\})$.
\item[(5)]If $\ell > 2n$, 
then a generic mapping $F_{\pi}\circ f:N\to \mathbb{R}^{\ell}$ is injective. 
\end{enumerate}
Moreover, by combining the assertions (3) and (5), 
for a given embedding $f:N\to U$, 
the following assertion (6) is obtained.
\begin{enumerate}
\item[(6)]If $\ell > 2n$ and $N$ is compact, 
then a generic mapping $F_{\pi}\circ f:N\to \mathbb{R}^{\ell}$ is an embedding. 
\end{enumerate}
\par 
\bigskip 
In Section \ref{section 2}, some standard definitions are reviewed, 
and the two main theorems (Theorems \ref{main} and \ref{main2}) are stated. 
Section \ref{section 3} (resp., Section \ref{section 4}) is devoted to
the proof of Theorem \ref{main} (resp., Theorem \ref{main2}). 
In Section \ref{section 5}, 
the assertions (1)-(6) above are shown. 
Moreover, in Section \ref{section 6}, as further applications, 
the two main theorems are adapted to quadratic mappings 
of $\mathbb{R}^m$ into $\mathbb{R}^\ell$ 
of a special type called ``generalized distance-squared mappings'' 
(for the precise definition of generalized distance-squared mappings, see Section \ref{section 6}). 
Since some corollaries in this paper 
(the assertion (6) in Section \ref{section 1}, Corollary \ref{embedding} 
in Section \ref{section 5} and 
Corollary \ref{stable} in Section \ref{section 6})
are also obtained by using the main theorem in \cite{F1}, which is an improvement of 
the main theorem in \cite{GP}, 
for the sake of readers' convenience, Section \ref{section 7} explains 
the main theorems in \cite{F1} and \cite{GP} as an appendix. 

\section{Preliminaries and the statements of Theorems \ref{main} and \ref{main2}}\label{section 2}

Let $N$ and $P$ be manifolds. Firstly, we recall the definition of transversality. 
\begin{definition}\label{transverse}
{\rm
Let $W$ be a submanifold of $P$. 
Let $g:N\to P$ be a mapping. 
\begin{enumerate}
\item 
We say that $g:N\to P$ is {\it transverse} to $W$ 
{\it at $q$} if $g(q)\not\in W$ or  
in the case of $g(q)\in W$, the following holds: 
\begin{eqnarray*}
dg_q(T_qN)+T_{g(q)}W=T_{g(q)}P.
\end{eqnarray*}
\item 
We say that $g:N\to P$ is {\it transverse} to $W$ 
if for any $q\in N$, 
the mapping $g$ is transverse to $W$ at $q$. 
\end{enumerate}
}
\end{definition}


We say that $g : N\to P$ is {\it $\mathcal{A}$-equivalent} to $h : N\to P$ 
if there exist diffeomorphisms $\Phi:N\to N$ and 
$\Psi:P\to P$ such that $g=\Psi \circ h \circ \Phi^{-1}$. 

Let $J^r(N,P)$ be the space of 
$r$-jets of mappings of $N$ into $P$. 
For a given mapping $g:N\to P$, 
the mapping $j^r g:N\to J^r(N,P)$ 
is defined by $q \mapsto j^r g(q)$  
(for details on the space $J^r(N,P)$ or 
the mapping $j^r g:N\to J^r(N,P)$, see for example, \cite{GG}). 

For the statement and the proof of Theorem \ref{main},
 it is sufficient to consider the case of $r=1$ and $P=\mathbb{R}^\ell$. 
Let $\{(U_\lambda ,\varphi _\lambda )\}_{\lambda \in \Lambda}$ be a coordinate neighborhood system of $N$. 
Let $\Pi :J^1(N,\mathbb{R}^\ell)$$\to N\times \mathbb{R}^\ell$ be the natural projection defined by $\Pi(j^1g(q))=(q,g(q))$. 
Let $\Phi _\lambda :\Pi^{-1}(U_\lambda \times \mathbb{R}^\ell)\to \varphi _\lambda (U_\lambda)\times \mathbb{R}^\ell \times J^1(n,\ell)$ be 
the homeomorphism defined by 
\begin{eqnarray*} 
\Phi _\lambda \left(j^1g(q)\right)=\left(\varphi _\lambda (q),g(q),j^1(\psi_{_\lambda} \circ g\circ \varphi _\lambda ^{-1}\circ \widetilde{\varphi} _\lambda)(0)\right), 
\end{eqnarray*} 
where 
$J^1(n, \ell)=\{ j^1g(0) \mid g : (\mathbb{R}^n, 0) \to  (\mathbb{R}^\ell, 0) \}$ 
and 
$\widetilde{\varphi} _\lambda : \mathbb{R}^n\to \mathbb{R}^n$ 
(resp., $\psi_{\lambda} : \mathbb{R}^m \to \mathbb{R}^m $) is the 
translation defined by 
$\widetilde{\varphi} _\lambda(0)=\varphi _\lambda (q)$ 
(resp., $\psi_{\lambda}(g(q))=0$). 
Then, $\{(\Pi^{-1}(U_\lambda \times \mathbb{R}^\ell), 
\Phi _\lambda )\}_{\lambda \in \Lambda}$ 
is a coordinate neighborhood system of $J^1(N,\mathbb{R}^\ell)$. 
A subset $X$ of $J^1(n,\ell)$ is said to be {\it $\mathcal{A}^1$-invariant} 
if for any $j^1g(0)\in X$, and for any two germs of diffeomorphisms  
$H : (\mathbb{R}^\ell, 0)\to (\mathbb{R}^\ell, 0)$ and 
$h : (\mathbb{R}^n, 0)\to (\mathbb{R}^n, 0)$, 
\<we have \> $j^1(H\circ g\circ h^{-1})(0)\in X$. 
Let $X$ be an $\mathcal{A}^1$-invariant submanifold of $J^1(n,\ell)$. 
Set 
\begin{eqnarray*}
X(N,\mathbb{R}^\ell)=\bigcup_{\lambda \in \Lambda}\Phi ^{-1}_\lambda \left(\varphi _\lambda (U_\lambda )\times \mathbb{R}^\ell \times X\right). 
\end{eqnarray*}
Then, the set $X(N,\mathbb{R}^\ell)$ is a subfiber-bundle of $J^1(N,\mathbb{R}^\ell)$ 
with the fiber $X$ such that 
\begin{eqnarray*}
{\rm codim}\ X(N,\mathbb{R}^\ell)&=&{\rm dim}\ J^1(N,\mathbb{R}^\ell)- 
{\rm dim}\ X(N,\mathbb{R}^\ell) \\
&=&{\rm dim}\ J^1(n,\ell)- {\rm dim}\ X \\
&=&{\rm codim}\ X.
\end{eqnarray*}

Then, the first main theorem in this paper is the following. 
\begin{theorem}\label{main}
Let $N$ be a manifold of dimension $n$. 
Let $f$ be an immersion 
of $N$ into an open subset $U$ of $\mathbb{R}^m$. 
Let $F:U\to \mathbb{R}^\ell$ be a mapping. 
If $X$ is an $\mathcal{A}^1$-invariant submanifold of $J^1(n,\ell)$, 
then there exists a subset $\Sigma$ with Lebesgue measure zero 
of $\mathcal{L}(\mathbb{R}^{m},\mathbb{R}^{\ell})$ 
such that for any $\pi \in \mathcal{L}(\mathbb{R}^{m},\mathbb{R}^{\ell})-\Sigma$, 
the mapping $j^1(F_\pi \circ f):
N\to J^1(N,\mathbb{R}^\ell)$ is transverse to the submanifold $X(N,\mathbb{R}^\ell)$.
\end{theorem}

Now, in order to state the second main theorem (Theorem \ref{main2}), 
we will prepare some definitions. 
Set $N^{(s)}=\{(q_1,\<q_2,\>\ldots ,q_s)\in N^s\mid q_i\not=q_j \ (i\not=j)\}$. 
Notice that $N^{(s)}$ is an open submanifold of $N^s$. 
For any mapping $g:N\to P$, 
let $g^{(s)}:N^{(s)}\to P^s$ be the mapping defined by 
\begin{eqnarray*}
g^{(s)}(q_1,\<q_2,\>\ldots ,q_s)=(g(q_1),\<g(q_2),\>\ldots ,g(q_s)). 
\end{eqnarray*}
Set 
$\Delta_s=\{(y,\ldots ,y) \in P^s \mid y\in P \}$. 
It is clearly seen that $\Delta_s$ is a submanifold of $P^s$ 
such that 
\begin{eqnarray*}
{\rm codim}\ \Delta_s ={\rm dim}\ P^s- 
{\rm dim}\ \Delta_s  =(s-1){\rm dim}\ P.
\end{eqnarray*}

\begin{definition}\label{normal crossings}
{\rm
Let $g$ be a mapping of $N$ into $P$. 
Then, $g$ is called a {\it  mapping with normal crossings} if  
for any positive integer $s$ $(s\geq 2)$, 
the mapping $g^{(s)}:N^{(s)}\to P^s$ is 
transverse to the submanifold $\Delta_s$. 
}
\end{definition}
For any injection $f:N\to \mathbb{R}^m$, set 
\begin{eqnarray*}
s_f={\rm max}\left\{s\ \middle| \ \forall (q_1,\<q_2,\>\ldots ,q_s)\in N^{(s)}, {\rm dim\ }\sum_{i=2}^s \mathbb{R}\overrightarrow{f(q_1)f(q_i)}=s-1 \right\}. 
\end{eqnarray*}
Since the mapping $f$ is injective, we get $2 \leq s_f$. 
Since $f(q_1),f(q_2),\ldots ,f(q_{s_f})$ are points of $\mathbb{R}^m$, 
it follows that $s_f\leq m+1$. 
Thus, we have  
\begin{eqnarray*}
2\leq s_f \leq m+1. 
\end{eqnarray*}
Furthermore, in the following, for a set $X$, 
we denote the number of its elements (or its cardinality) by $|X|$. 
Then, the second main theorem in this paper is the following. 
\begin{theorem}\label{main2}
Let $N$ be a manifold of dimension $n$. 
Let $f$ be an injection of $N$ into an open subset 
$U$ of $\mathbb{R}^m$. 
Let $F:U\to \mathbb{R}^\ell$ be a mapping. 
Then, there exists a subset $\Sigma$ of $\mathcal{L}(\mathbb{R}^m, \mathbb{R}^\ell)$ with Lebesgue measure zero such that 
for any $\pi \in \mathcal{L}(\mathbb{R}^m, \mathbb{R}^\ell)-\Sigma $, \<and \> 
for any $s$ $(2\leq s \leq s_f)$, 
the mapping $(F_\pi \circ f)^{(s)}:N^{(s)}\to (\R^\ell)^s$ is transverse to the 
submanifold $\Delta_s$. 
Moreover, if the mapping $F_\pi$ 
satisfies that $| F_\pi^{-1}(y) | \leq s_f$ for any $y\in\mathbb{R}^\ell$, 
then $F_\pi \circ f:N\to \mathbb{R}^\ell$ is a mapping with normal crossings.    
\end{theorem}

The following well known result is important for the proofs of Theorems \ref{main} and \ref{main2}. 
\begin{lemma}[\cite{abra}, \cite{GP}]\label{abra}
Let $N$, $P$, $Z$ be manifolds, and 
let $W$ be a submanifold of $P$. Let $\Gamma:N\times Z\to P$ be a mapping.
If $\Gamma$ is transverse to $W$, then there exists a subset $\Sigma$ of $Z$ with 
Lebesgue measure zero such that for any $p \in Z-\Sigma $, 
the mapping 
$\Gamma_p:N\to P$ is transverse to $W$, 
where $\Gamma_p(q)=\Gamma(q,p)$. 
\end{lemma}

\begin{remark}
{\rm 
\begin{enumerate}
\item 
We explain the advantage that the domain of the mapping $F$ 
is an arbitrary open set. 
Suppose that $U=\mathbb{R}$. 
Let $F:\mathbb{R}\to \mathbb{R}$ be the mapping defined by $x\mapsto |x|$. 
Since $F$ is not differentiable at $x=0$, we cannot apply Theorems \ref{main} and \ref{main2} to 
the mapping $F:\mathbb{R}\to \mathbb{R}$. 

On the other hand, if $U=\mathbb{R}-\{0\}$, then 
Theorems \ref{main} and \ref{main2} 
can be applied to the restriction $F|_U$.
\item 
There is a case of $s_f =3$ as follows. 
If $n+1 \leq m$, $N=S^n$ and 
$f:S^n\to \mathbb{R}^m$ is the inclusion $f(x)=(x,0,\ldots ,0)$, 
then it is easily seen that $s_f=3$. 
Indeed, suppose that there exists a point $(q_1, q_2, q_3)\in (S^n)^{(3)}$ 
such that ${\rm dim\ }\sum_{i=2}^3 \mathbb{R}\overrightarrow{f(q_1)f(q_i)}=1$. 
Then, since the number of the intersections of $f(S^n)$ and 
a straight line of $\mathbb{R}^m$ is at most two, 
this contradicts the assumption. 
Thus, 
we get $s_f\geq3$. 
From $S^1\times \{0\} \subset f(S^n)$, 
it follows that $s_f<4$, where $0=\underbrace{(0,\ldots ,0)}_{(m-2)
\text{-tuple}}$. 
Hence, we have $s_f=3$.   
\item 

The essential idea for the proofs of Theorems \ref{main} and \ref{main2}  
is to apply Lemma \ref{abra}, and it 
is almost similar to the idea of the proofs of main results in \cite{G3}. 
Nevertheless, the two main theorems in this paper are drastically improved. 
As an effect of the improvement, 
many applications 
are obtained by the two main theorems 
(for the applications, see Sections \ref{section 5} and \ref{section 6}). 
\end{enumerate}
}
\end{remark}

\section{Proof of Theorem \ref{main}}\label{section 3}

Let $(\alpha_{ij})_{1\leq i \leq \ell, 1\leq j \leq m}$ be a representing matrix of 
a linear mapping $\pi:\mathbb{R}^m\to \mathbb{R}^\ell$. 
Set $F_{\alpha}=F_{\pi}$, and we have
\[
F_{\alpha}(x)=
\biggl(F_{1}(x)+\sum_{j=1}^{m}\alpha_{1j}x_j, 
F_{2}(x)+\sum_{j=1}^{m}\alpha_{2j}x_j, 
\ldots,
F_{\ell}(x)+\sum_{j=1}^{m}\alpha_{\ell j}x_j\biggr), 
\eqno (3.1)
\]
where $F=(F_1, F_2,\ldots ,F_\ell)$, 
 $\alpha=(\alpha_{11},\<\alpha_{12},\ldots, 
 \alpha_{1m},\ldots ,
 \alpha_{\ell 1},\alpha_{\ell 2},\ldots ,
\alpha_{\ell m})\in (\mathbb{R}^m)^\ell$ and $x=(x_1,x_2,\ldots ,x_m)$. 
For a given immersion $f:N\to U$, 
the mapping $F_{\alpha}\circ f:N\to \mathbb{R}^\ell$ is 
given as follows: 
\[
F_{\alpha}\circ f=
\biggl(F_{1}\circ f+\sum_{j=1}^{m}\alpha_{1j}f_j, 
F_{2}\circ f+\sum_{j=1}^{m}\alpha_{2j}f_j, 
\ldots,
F_{\ell}\circ f+\sum_{j=1}^{m}\alpha_{\ell j}f_j\biggr),
\eqno (3.2)
\]
where $f=(f_1,f_2,\ldots ,f_m)$.  
Since we have the natural identification 
$\mathcal{L}(\mathbb{R}^{m},\mathbb{R}^{\ell})=(\mathbb{R}^{m})^\ell$, 
in order to prove Theorem \ref{main}, 
it is sufficient to show that 
there exists a subset $\Sigma$ 
with Lebesgue measure zero of $(\mathbb{R}^m)^{\ell}$ 
such that for any $\alpha \in (\mathbb{R}^m)^{\ell}-\Sigma$, 
the mapping $j^{1}(F_{\alpha}\circ f):
N\to J^{1}(N,\mathbb{R}^\ell)$ is transverse to 
the given submanifold $X(N,\mathbb{R}^\ell)$. 

Now, let $\Gamma:N\times (\mathbb{R}^m)^\ell \to J^1(N,\mathbb{R}^\ell)$ 
be the mapping defined by 
\[
\Gamma(q,\alpha)=j^1(F_\alpha \circ f)(q).
\]
If the mapping $\Gamma$ is transverse to the submanifold 
$X(N,\mathbb{R}^\ell)$, 
then \<from \> Lemma \ref{abra},  
it follows that there exists a subset $\Sigma$ of $(\mathbb{R}^m)^\ell$ with 
Lebesgue measure zero such that 
for any $\alpha \in (\mathbb{R}^m)^\ell-\Sigma $, 
the mapping $\Gamma_\alpha :N\to J^1(N,\mathbb{R}^\ell)$ 
($\Gamma_\alpha=j^1(F_\alpha \circ f)$) 
is transverse to the submanifold 
$X(N,\mathbb{R}^\ell)$. 
Thus, in order to finish the proof of Theorem \ref{main}, 
it is sufficient to show that if $\Gamma(\widetilde{q},\widetilde{\alpha})
\in X(N,\mathbb{R}^\ell)$, then the following holds:
\[
d\Gamma_{(\widetilde{q},\widetilde{\alpha})}(T_{(\widetilde{q},\widetilde{\alpha})}
(N\times (\mathbb{R}^m)^\ell))
+
T_{\Gamma(\widetilde{q},\widetilde{\alpha})}X(N,\mathbb{R}^\ell)
=
T_{\Gamma(\widetilde{q},\widetilde{\alpha})}J^1(N,\mathbb{R}^\ell). \eqno (3.3)
\]
As in Section \ref{section 2}, 
let $\{(U_\lambda ,\varphi _\lambda )\}_{\lambda \in \Lambda}$ 
(resp., $\{(\Pi^{-1}(U_\lambda \times \mathbb{R}^\ell), 
\Phi _\lambda )\}_{\lambda \in \Lambda}$) 
be a coordinate neighborhood system of $N$ (resp., $J^1(N,\mathbb{R}^\ell)$). 
There exists a coordinate neighborhood 
$\left(U_{\widetilde{\lambda}}\times (\mathbb{R}^m)^\ell, \varphi_{\widetilde{\lambda}}\times id \right)$ 
containing the point $(\widetilde{q},\widetilde{\alpha})$ of 
$N\times (\mathbb{R}^m)^\ell$, 
where $id$ is the identity mapping of $(\mathbb{R}^m)^\ell$ into $(\mathbb{R}^m)^\ell$, 
and the mapping $\varphi_{\widetilde{\lambda}}\times id : 
U_{\widetilde{\lambda}}\times (\mathbb{R}^m)^\ell \to 
\varphi_{\widetilde{\lambda}}(U_{\widetilde{\lambda}})\times (\mathbb{R}^m)^\ell$ ($\subset \mathbb{R}^n\times (\mathbb{R}^m)^\ell$) 
is defined by 
$\left(\varphi_{\widetilde{\lambda}}\times id\right)(q,\alpha)=
\left(\varphi_{\widetilde{\lambda}}(q), id(\alpha)\right)$.
There exists a coordinate neighborhood 
$\left(\Pi^{-1}(U_{\widetilde{\lambda}} \times \mathbb{R}^\ell), 
\Phi _{\widetilde{\lambda}} \right)$ 
containing the point $\Gamma (\widetilde{q},\widetilde{\alpha})$ of 
$J^1(N,\mathbb{R}^\ell)$. 
Let $t=(t_1,\<t_2,\>\ldots ,t_n)\in \mathbb{R}^n$ be a local coordinate on 
$\varphi_{\widetilde{\lambda}}(U_{\widetilde{\lambda}})$ 
containing 
$\varphi_{\widetilde{\lambda}}(\widetilde{q})$. 
Then, the mapping $\Gamma$ is locally given by the following:
{\small 
\begin{eqnarray*}
&{}&(\Phi _{\widetilde{\lambda}} \circ \Gamma \circ 
(\varphi_{\widetilde{\lambda}}\times id )^{-1})(t,\alpha)\\
&=&(\Phi _{\widetilde{\lambda}}\circ j^1(F_\alpha \circ f)\circ \varphi_{\widetilde{\lambda}}^{-1})(t)\\
&=&\left(t, (F_\alpha \circ f \circ \varphi_{\widetilde{\lambda}}^{-1})(t), \right.  \\
&& \frac{\partial (F_{\alpha, 1}\circ f\circ \varphi_{\widetilde{\lambda}}^{-1})}
{\partial t_1}(t), 
\<
\frac{\partial (F_{\alpha, 1}\circ f\circ \varphi_{\widetilde{\lambda}}^{-1})}
{\partial t_2}(t), 
\>
\ldots ,
\frac{\partial (F_{\alpha, 1}\circ f\circ \varphi_{\widetilde{\lambda}}^{-1})}
{\partial t_n}(t),
\\ 
&& \< 
\frac{\partial (F_{\alpha, 2}\circ f\circ \varphi_{\widetilde{\lambda}}^{-1})}
{\partial t_1}(t), 
\frac{\partial (F_{\alpha, 2}\circ f\circ \varphi_{\widetilde{\lambda}}^{-1})}
{\partial t_2}(t), 
\ldots ,
\frac{\partial (F_{\alpha, 2}\circ f\circ \varphi_{\widetilde{\lambda}}^{-1})}
{\partial t_n}(t), 
\> 
\\ 
\\
&&\hspace{120pt}\cdots \cdots \cdots , \\ 
\\
&& \left. \frac{\partial (F_{\alpha, \ell}\circ f\circ \varphi_{\widetilde{\lambda}}^{-1})}
{\partial t_1}(t)
,
\<
\frac{\partial (F_{\alpha, \ell}\circ f\circ \varphi_{\widetilde{\lambda}}^{-1})}
{\partial t_2}(t), 
\>
\ldots ,
\frac{\partial (F_{\alpha, \ell}\circ f\circ \varphi_{\widetilde{\lambda}}^{-1})}
{\partial t_n}(t)
\right)\\
&=&\left(t, (F_\alpha \circ f \circ \varphi_{\widetilde{\lambda}}^{-1})(t), \right. \\
&&\frac{\partial F_1\circ \widetilde{f}}{\partial t_1}(t)+
\sum_{j=1}^m \alpha_{1j}\frac{\partial \widetilde{f}_j}{\partial t_1}(t),
\<
\frac{\partial F_1\circ \widetilde{f}}{\partial t_2}(t)+
\sum_{j=1}^m \alpha_{1j}\frac{\partial \widetilde{f}_j}{\partial t_2}(t),
\>
\ldots ,
\frac{\partial F_1\circ \widetilde{f}}{\partial t_n}(t)+
\sum_{j=1}^m \alpha_{1j}\frac{\partial \widetilde{f}_j}{\partial t_n}(t),
\\ 
&& 
\<
\frac{\partial F_2\circ \widetilde{f}}{\partial t_1}(t)+
\sum_{j=1}^m \alpha_{2j}\frac{\partial \widetilde{f}_j}{\partial t_1}(t),
\frac{\partial F_2\circ \widetilde{f}}{\partial t_2}(t)+
\sum_{j=1}^m \alpha_{2j}\frac{\partial \widetilde{f}_j}{\partial t_2}(t),
\ldots ,
\frac{\partial F_2\circ \widetilde{f}}{\partial t_n}(t)+
\sum_{j=1}^m \alpha_{2j}\frac{\partial \widetilde{f}_j}{\partial t_n}(t),
\> 
\\
\\
&&\hspace{160pt}\cdots \cdots \cdots , 
\\
\\
&& \left. 
\frac{\partial F_\ell \circ \widetilde{f}}{\partial t_1}(t)+\sum_{j=1}^m \alpha_{\ell j}\frac{\partial \widetilde{f}_j}{\partial t_1}(t),
\<
\frac{\partial F_\ell \circ \widetilde{f}}{\partial t_2}(t)+\sum_{j=1}^m \alpha_{\ell j}\frac{\partial \widetilde{f}_j}{\partial t_2}(t),
\>
\ldots ,
\frac{\partial F_\ell \circ \widetilde{f}}{\partial t_n}(t)+
\sum_{j=1}^m \alpha_{\ell j}\frac{\partial \widetilde{f}_j}{\partial t_n}(t)\right), 
\end{eqnarray*} 
}where $F_\alpha=(F_{\alpha, 1},\<F_{\alpha, 2},\> \ldots , F_{ \alpha, \ell})$ and 
$\widetilde{f}=(\widetilde{f}_1, \<\widetilde{f}_2,\> \ldots ,\widetilde{f}_m)
=(f_1\circ \varphi_{\widetilde{\lambda}}^{-1},
\<
f_2\circ \varphi_{\widetilde{\lambda}}^{-1},
\>
\ldots , f_m\circ
\varphi_{\widetilde{\lambda}}^{-1})=f\circ \varphi_{\widetilde{\lambda}}^{-1}$. 
The Jacobian matrix of the mapping $\Gamma$ at 
$(\widetilde{q},\widetilde{\alpha})$ is the following:\\
\begin{eqnarray*}
J\Gamma_{(\widetilde{q},\widetilde{\alpha})}=
\left(
\begin{array}{@{\,}c@{\,\,}|@{\,\,}c@{\,\,\,}c@{\,\,\,}c@{\,\,\,}c@{\,\,\,}c}
E_n       &    0    &         \cdots          &   \cdots        &0    \\
\hline                &   \ast      &    \cdots    &  \cdots      & \ast  \\ 
      &   {}^t\!(Jf_{\widetilde{q}}) & & \bigzerol &\\
 \ast    &              & {}^t\!(Jf_{\widetilde{q}})    &       &  \\
    &    & \bigzerol &\ddots &      \\
      &    &  &  &   {}^t\!(Jf_{\widetilde{q}}) \\
\end{array}
\right)_{(t, \alpha)=(\varphi_{\widetilde{\lambda}}(\widetilde{q}),\widetilde{\alpha})}, 
\end{eqnarray*}
where $E_n$ is the $n\times n$ unit matrix and 
$Jf_{\widetilde{q}}$ is the Jacobian matrix of the mapping $f$ at $\widetilde{q}$. 
Note that ${}^t\!(Jf_{\widetilde{q}})$ is the transpose of the matrix 
$Jf_{\widetilde{q}}$ and that there are $\ell$ copies of ${}^t\!(Jf_{\widetilde{q}})$ 
in the above description of $J\Gamma_{(\widetilde{q},\widetilde{\alpha})}$. 
Since $X(N,\mathbb{R}^\ell)$ is a subfiber-bundle of $J^1(N,\mathbb{R}^\ell)$ 
with the fiber $X$, 
it is clear that in order to show (3.3), it 
suffices to prove that the matrix $M_1$ given below 
has rank $n+\ell+n\ell$: 
\begin{eqnarray*}
M_1=
\left(
\begin{array}{@{\,}c@{\,\,}|@{\,\,}c@{\,\,\,}c@{\,\,\,}c@{\,\,\,}c@{\,\,\,}c}
E_{n+\ell}       &    \ast    &         \cdots          &   \cdots        &\ast   \\
\hline     \\[-3.3mm]   &   {}^t\!(Jf_{\widetilde{q}})  & & \bigzerol &\\
\bigzerol  &              & {}^t\!(Jf_{\widetilde{q}})      &       &  \\
    &    & \bigzerol &\ddots &      \\
      &    &  &  &   {}^t\!(Jf_{\widetilde{q}})  \\
\end{array}
\right)_{(t, \alpha)=(\varphi_{\widetilde{\lambda}}(\widetilde{q}),\widetilde{\alpha})}, 
\end{eqnarray*}
where $E_{n+\ell}$ is the $(n+\ell)\times (n+\ell)$ unit matrix. 
Note that there are $\ell$ copies of ${}^t\!(Jf_{\widetilde{q}})$ 
in the above description of $M_1$. 
Notice that for any $i$ $(1\leq i \leq m\ell)$, the $(n+\ell+i)$-th 
column vector of $M_1$ coincides with 
the $(n+i)$-th column vector of $J\Gamma_{(\widetilde{q},\widetilde{\alpha})}$. 
Since 
the mapping $f$ is an immersion $(n\leq m)$, we have that the rank of the matrix 
$M_1$ 
is equal to $n+\ell+n\ell$. 
Hence, we have (3.3). 
\hfill\qed
\section{Proof of Theorem \ref{main2}}\label{section 4}
By the same method as in the proof of Theorem \ref{main}, 
set $F_{\alpha}=F_{\pi}$, where $F_{\alpha}$ is given by (3.1) 
in Section \ref{section 3}. 
For a given injection $f:N\to U$, 
the mapping $F_{\alpha}\circ f:N\to \mathbb{R}^\ell$ is 
given by the same expression as (3.2).  
Since we have the natural identification 
$\mathcal{L}(\mathbb{R}^{m},\mathbb{R}^{\ell})=(\mathbb{R}^{m})^\ell$, 
in order to show that 
there exists a subset $\Sigma$ of $\mathcal{L}(\mathbb{R}^m, \mathbb{R}^\ell)$ with Lebesgue measure zero such that 
for any $\pi \in \mathcal{L}(\mathbb{R}^m, \mathbb{R}^\ell)-\Sigma $, 
\<and \> 
for any $s$ $(2\leq s \leq s_f)$, 
the mapping $(F_\pi \circ f)^{(s)}:N^{(s)}\to (\R^\ell)^s$ is transverse to the 
submanifold $\Delta_s$,  
it is sufficient to show that 
there exists a subset $\Sigma$ of $(\mathbb{R}^m)^{\ell}$ 
with Lebesgue measure zero 
such that for any $\alpha \in (\mathbb{R}^m)^{\ell}-\Sigma$, 
\<and \> 
for any $s$ $(2\leq s \leq s_f)$, 
the mapping $(F_{\alpha}\circ f)^{(s)}:
N^{(s)}\to (\mathbb{R}^\ell)^s$ is transverse to $\Delta_s$.  

Now, let $s$ be a positive integer satisfying $2\leq s \leq s_f$.  
Let $\Gamma : N^{(s)}\times (\mathbb{R}^m)^\ell \to (\mathbb{R}^{\ell})^s$ be the mapping defined by 
\begin{eqnarray*}
\Gamma(q_1,\<q_2,\>\ldots ,q_s,\alpha)=
\left( (F_\alpha \circ f)(q_1), \< (F_\alpha \circ f)(q_2),\>\ldots, 
(F_\alpha \circ f)(q_s)\right).
\end{eqnarray*}
If for any positive integer $s$ $(2\leq s \leq s_f)$, 
the mapping $\Gamma$ is transverse to $\Delta_s$, then \<from \> Lemma \ref{abra}, 
 it follows that for any positive integer $s$ $(2\leq s \leq s_f)$, 
 there exists a subset $\Sigma_s$ of $(\mathbb{R}^m)^\ell$ 
 with Lebesgue measure zero 
such that for any $\alpha  \in (\mathbb{R}^m)^{\ell}-\Sigma_s$,  
the mapping $\Gamma_\alpha : 
N^{(s)}\to (\mathbb{R}^\ell)^s$ $(\Gamma_\alpha=(F_\alpha\circ f)^{(s)})$ 
is transverse to $\Delta_s$. 
Then, set $\Sigma=\bigcup_{s=2}^{s_f}\Sigma_s$. 
It is clearly seen that $\Sigma$ is a subset of  $(\mathbb{R}^m)^\ell$ 
with Lebesgue measure zero. 
Therefore, it follows that  
for any $\alpha \in (\mathbb{R}^m)^{\ell}-\Sigma $, and  
for any $s$ $(2\leq s \leq s_f)$, 
the mapping $\Gamma_\alpha : 
N^{(s)}\to (\mathbb{R}^\ell)^s$ $(\Gamma_\alpha=(F_\alpha\circ f)^{(s)})$ 
is transverse to $\Delta_s$. 

Hence, for the proof, 
it is sufficient to show that 
for any positive integer $s$ $(2\leq s \leq s_f)$,   
if $\Gamma(\widetilde{q}, \widetilde{\alpha})\in \Delta_s $ 
$(\widetilde{q}=(\widetilde{q}_1,
\widetilde{q}_2,
\ldots ,\widetilde{q}_s))$, 
then the following holds: 
\[
d\Gamma_{(\widetilde{q}, \widetilde{\alpha})}(T_{(\widetilde{q}, \widetilde{\alpha})}(N^{(s)}\times (\mathbb{R}^m)^\ell))+
T_{\Gamma(\widetilde{q}, \widetilde{\alpha})}\Delta_s
=T_{\Gamma(\widetilde{q}, \widetilde{\alpha})}(\mathbb{R}^\ell)^s. \eqno (4.1)
\]
Let $\{(U_\lambda ,\varphi _\lambda )\}_{\lambda \in \Lambda}$ be a coordinate  neighborhood system of $N$. 
There exists a coordinate neighborhood 
$(U_{\widetilde{\lambda}_1}\times \<U_{\widetilde{\lambda}_2}\times \>
\cdots \times U_{\widetilde{\lambda}_s} \times (\mathbb{R}^m)^\ell, 
\varphi_{\widetilde{\lambda}_1}\times 
\<\varphi_{\widetilde{\lambda}_2}\times 
\>
\cdots \times \varphi_{\widetilde{\lambda}_s}\times id )$ 
containing the point $(\widetilde{q}, \widetilde{\alpha})$ 
of $N^{(s)}\times (\mathbb{R}^m)^\ell$, 
where $id$ is the identity mapping of $(\mathbb{R}^m)^\ell$ into $(\mathbb{R}^m)^\ell$, 
and the mapping 
$\varphi_{\widetilde{\lambda}_1}\times 
\<
\varphi_{\widetilde{\lambda}_2}\times 
\>
\cdots \times \varphi_{\widetilde{\lambda}_s}\times id  : 
U_{\widetilde{\lambda}_1}\times 
\<
U_{\widetilde{\lambda}_2}\times 
\>
\cdots \times U_{\widetilde{\lambda}_s} \times (\mathbb{R}^m)^\ell 
 \to (\mathbb{R}^n)^s\times (\mathbb{R}^m)^\ell$ is defined by 
$(\varphi_{\widetilde{\lambda}_1}\times 
\<
\varphi_{\widetilde{\lambda}_2}\times 
\>
\cdots \times \varphi_{\widetilde{\lambda}_s}\times id )(q_1,\<q_2,\>\ldots ,q_s,\alpha)=
(\varphi_{\widetilde{\lambda}_1}(q_1),  
\<
\varphi_{\widetilde{\lambda}_2}(q_2),  
\>
\ldots ,
\varphi_{\widetilde{\lambda}_s}(q_s), id(\alpha ))$.
Let $t_i=(t_{i1},t_{i2},\ldots ,t_{in})$ be a local coordinate 
around $\varphi_{\widetilde{\lambda}_i}(\widetilde{q}_i)$ $(1\leq i \leq s)$.  
Then, the mapping $\Gamma$ is locally given by the following:
\begin{eqnarray*}
&&\Gamma \circ \left(\varphi_{\widetilde{\lambda}_1}\times 
\<
\varphi_{\widetilde{\lambda}_2}\times 
\>
\cdots \times \varphi_{\widetilde{\lambda}_s}\times id\right)^{-1}
(t_1,\<t_2,\>\ldots ,t_s,\alpha)
\\
&=&\left( (F_\alpha \circ f\circ \varphi_{\widetilde{\lambda}_1}^{-1})(t_1), 
\<
(F_\alpha \circ f\circ \varphi_{\widetilde{\lambda}_2}^{-1})(t_2), 
\>
\ldots ,
(F_\alpha \circ f\circ \varphi_{\widetilde{\lambda}_s}^{-1})(t_s) 
\right)
\\
&=&\left(
F_1\circ \widetilde{f}(t_1)+\sum_{j=1}^m\alpha_{1j}\widetilde{f}_j(t_1), 
\<
F_2\circ \widetilde{f}(t_1)+\sum_{j=1}^m\alpha_{2j}\widetilde{f}_j(t_1), 
\>
\ldots ,
F_\ell \circ \widetilde{f}(t_1)+\sum_{j=1}^m\alpha_{\ell j}\widetilde{f}_j(t_1), \right.\\
&& 
\<
F_1\circ \widetilde{f}(t_2)+\sum_{j=1}^m\alpha_{1j}\widetilde{f}_j(t_2), 
F_2\circ \widetilde{f}(t_2)+\sum_{j=1}^m\alpha_{2j}\widetilde{f}_j(t_2), 
\ldots ,
F_\ell \circ \widetilde{f}(t_2)+\sum_{j=1}^m\alpha_{\ell j}\widetilde{f}_j(t_2), 
\>  
\\
\\
&&\hspace{150pt}\cdots \cdots \cdots , 
\\
\\
&& \left. 
F_1\circ \widetilde{f}(t_s)+\sum_{j=1}^m\alpha_{1j}\widetilde{f}_j(t_s),
\<
F_2\circ \widetilde{f}(t_s)+\sum_{j=1}^m\alpha_{2j}\widetilde{f}_j(t_s),
\>
\ldots ,
F_\ell \circ \widetilde{f}(t_s)+\sum_{j=1}^m\alpha_{\ell j}\widetilde{f}_j(t_s)
\right), 
\end{eqnarray*}  
where
$\widetilde{f}(t_i)=(\widetilde{f}_1(t_i),
\<
\widetilde{f}_2(t_i),
\>
\ldots ,\widetilde{f}_m(t_i))
=(f_1\circ \varphi_{\widetilde{\lambda}_i}^{-1}(t_i),
\<
f_2\circ \varphi_{\widetilde{\lambda}_i}^{-1}(t_i),
\>
\ldots ,
f_m\circ \varphi_{\widetilde{\lambda}_i}^{-1}(t_i))$ $(1\leq i\leq s)$. 
For simplicity, set $t=(t_1,\<t_2,\>\ldots ,t_s)$ and 
$z=(\varphi_{\widetilde{\lambda}_1}\times 
\<
\varphi_{\widetilde{\lambda}_2}\times 
\>
\cdots \times \varphi_{\widetilde{\lambda}_s})
(\widetilde{q}_1,\<\widetilde{q}_2,\>\ldots ,\widetilde{q}_s)$. 

The Jacobian matrix of the mapping $\Gamma$ at 
$(\widetilde{q}, \widetilde{\alpha})$ is the following: 
\begin{eqnarray*}
J\Gamma_{(\widetilde{q}, \widetilde{\alpha})}=
\left(
\begin{array}{@{\,\,\,}c@{\,\,\,\,\,}|@{\,\,\,\,\,}c@{\,\,\,}}
\ast & B(t_1) \\
\ast & B(t_2) \\
 \vdots & \vdots    \\ 
\ast & B(t_s) \\
\end{array}
\right)_{(t,\alpha)=(z, \widetilde{\alpha})},
\end{eqnarray*}
where 
\begin{eqnarray*}
B(t_i)=
\left. 
\left(
\begin{array}{ccccccc}
{\bf b}_{}(t_i)& & &\bigzerol \\
&{\bf b}_{}(t_i)&& \\
\bigzerol &   &  \ddots   &    \\ 
&&   &         &{\bf b}_{}(t_i)
\end{array}
\right)
\right\}
\,\text{{\rm $\ell$ rows}}
\end{eqnarray*}
and 
${\bf b}_{ }(t_i)=(
\widetilde{f}_1(t_i),
\widetilde{f}_2(t_i),
\ldots ,\widetilde{f}_m(t_i))$. 
By the construction of 
$T_{\Gamma(\widetilde{q}, \widetilde{\alpha})}\Delta_s$, 
in order to show (4.1), it is sufficient to 
show that the rank of the following matrix $M_2$ is equal to $\ell s$: 
\begin{eqnarray*}
M_2=
\left(
\begin{array}{@{\,\,\,}c@{\,\,\,\,\,}|@{\,\,\,\,\,}c@{\,\,\,}}
E_\ell & B(t_1) \\
E_\ell & B(t_2) \\
 \vdots & \vdots    \\ 
E_\ell & B(t_s) \\
\end{array}
\right)_{t=z}.
\end{eqnarray*}
There exists an 
$\ell s\times \ell s$ regular matrix $Q_1$ 
such that 
\begin{eqnarray*}
Q_1M_2=
\left(
\begin{array}{@{\,\,\,}c@{\,\,\,\,\,}|@{\,\,\,\,\,}c@{\,\,\,}}
E_\ell & B(t_1) \\
0 & B(t_2)- B(t_1) \\
 \vdots & \vdots    \\ 
0 & B(t_s)- B(t_1) \\
\end{array}
\right)_{t=z}.
\end{eqnarray*}
There exists an 
$(\ell+m\ell)\times (\ell+m\ell)$ regular matrix $Q_2$ 
such that 
\begin{eqnarray*}
Q_1M_2Q_2&=&
\left(
\begin{array}{@{\,\,\,}c@{\,\,\,\,}|@{\,\,\,\,}c@{\,\,\,}}
E_\ell & 0 \\
0 & B(t_2)- B(t_1) \\
 \vdots & \vdots    \\ 
0 & B(t_s)- B(t_1) \\
\end{array}
\right)_{t=z}
\\
&=&
\begin{array}{c@{\,\,\,\,\,}c@{\,\,\,\,}|@{\,\,\,\,}c@{\,\,\,\,}@{\,\,\,\,}c@{\,\,\,\,}c
@{\,\,\,\,}c@{\,\,\,\,}c@{\,\,\,\,}ccc}
\ldelim({17}{4pt}[] && & &  &&\rdelim){17}{4pt}[]\\
&E_\ell &&\bigzerol & & \\
&&     &    &    \\ 
 \cline{2-6}  &&&&&&&&\rdelim\}{6}{10pt}[$\ell$ rows]\\
\\[-7mm]
&&\overrightarrow{\widetilde{f}(t_1)\widetilde{f}(t_2)} & & &  \bigzerol  & && \\   
&\bigzerol&&\overrightarrow{\widetilde{f}(t_1)\widetilde{f}(t_2)}  & \\
&&\bigzerol& &\ddots   &    \\ 
&&  & &&  \overrightarrow{\widetilde{f}(t_1)\widetilde{f}(t_2)}   \\
&&  & && \\
\\[-8mm]
 \cline{2-6} 
&\vdots& \vdots &\vdots &\vdots& \vdots  \\
  \cline{2-6}  &&&&&&&&\rdelim\}{6}{10pt}[$\ell$ rows]\\
\\[-7mm]
 &&  \overrightarrow{\widetilde{f}(t_1)\widetilde{f}(t_s)}   &  & &   \bigzerol       \\
&\bigzerol&& \overrightarrow{\widetilde{f}(t_1)\widetilde{f}(t_s)}   & &\\
&&\bigzerol&   &  \ddots   &    \\ 
&&  & &&   \overrightarrow{\widetilde{f}(t_1)\widetilde{f}(t_s)}    \\
\end{array}, 
\end{eqnarray*}
where 
$\overrightarrow{\widetilde{f}(t_1)\widetilde{f}(t_i)} =(
\widetilde{f}_1(t_i)-\widetilde{f}_1(t_1), 
\widetilde{f}_2(t_i)-\widetilde{f}_2(t_1), 
\ldots ,\widetilde{f}_m(t_i)-\widetilde{f}_m(t_1) 
)$ $(2\leq i \leq s)$ and $t=z$. 
From  $s-1\leq s_f-1$ and the definition of $s_f$, it follows that 
\begin{eqnarray*}
{\rm dim} \sum_{i=2}^s\mathbb{R}\overrightarrow{\widetilde{f}(t_1)\widetilde{f}(t_i)}=s-1,
\end{eqnarray*}
where $t=z$. Thus, 
by the construction of the matrix 
$Q_1M_2Q_2$ and $s-1\leq m$, 
we have that the rank of the matrix 
$Q_1M_2Q_2$ 
is equal to $\ell s$. 
Hence, the rank of the matrix $M_2$ 
must be equal to $\ell s$. 
Therefore, we have (4.1). 
Thus, there exists a subset $\Sigma$ of $\mathcal{L}(\mathbb{R}^m, \mathbb{R}^\ell)$ with Lebesgue measure zero such that 
for any $\pi \in \mathcal{L}(\mathbb{R}^m, \mathbb{R}^\ell)-\Sigma $, 
\<and \> 
for any $s$ $(2\leq s \leq s_f)$, 
the mapping $(F_\pi \circ f)^{(s)}:N^{(s)}\to (\R^\ell)^s$ is transverse to the 
submanifold $\Delta_s$. 

Moreover, suppose that 
the mapping $F_\pi$ 
satisfies that $| F_\pi^{-1}(y) | \leq s_f$ for any $y\in\mathbb{R}^\ell$. 
Since $f:N\to \mathbb{R}^m$ is injective, 
it follows that $| (F_\pi \circ f)^{-1}(y) | \leq s_f$ for any $y\in\mathbb{R}^\ell$. 
Hence, it follows that for any positive integer $s$ with $s\geq s_f+1$, we have 
$(F_\pi \circ f)^{(s)}(N^{(s)})\bigcap \Delta_s=\emptyset$.  
Namely, for any positive integer $s$ with $s\geq s_f+1$, 
the mapping $(F_\pi \circ f)^{(s)}$ is transverse to $\Delta_s$. 
Thus, $F_\pi \circ f:N\to \mathbb{R}^\ell$ is a mapping with normal crossings. 
\hfill\qed
\section{Applications of Theorems \ref{main} and \ref{main2}}\label{section 5}
In Subsection \ref{application1} (resp., Subsection \ref{application2}), 
applications of Theorem \ref{main} (resp., Theorem \ref{main2}) 
are stated and proved. 
In Subsection \ref{application2}, applications obtained by combining Theorems \ref{main} and  \ref{main2}  
are also given. 
\subsection{Applications of Theorem \ref{main}}\label{application1}

Set 
\begin{eqnarray*} 
\Sigma ^k=\left\{j^1g(0)\in J^1(n,\ell)\mid {\rm corank\ }Jg(0)=k\right\}, 
\end{eqnarray*}
where ${\rm corank\ }Jg(0)={\rm min} \{n,\ell \}-{\rm rank\ } Jg(0)$ 
and $k=1,\<2,\>\ldots ,{\rm min}\{n, \ell\}$. 
Then, $\Sigma ^k$ is an $\mathcal{A}^1$-invariant submanifold of $J^1(n,\ell)$. 
Set 
\begin{eqnarray*}
\Sigma^k(N,\mathbb{R}^\ell)=\bigcup_{\lambda \in \Lambda}\Phi ^{-1}_\lambda \left(\varphi _\lambda (U_\lambda )\times \mathbb{R}^\ell \times \Sigma ^k\right), 
\end{eqnarray*}
where the mappings $\Phi _\lambda$ and $\varphi _\lambda $ are as defined 
in Section \ref{section 2}.  
Then, the set $\Sigma ^k(N,\mathbb{R}^\ell)$ is a subfiber-bundle of 
$J^1(N,\mathbb{R}^\ell)$ with the fiber $\Sigma^k$ such that 
\begin{eqnarray*}
{\rm codim}\ \Sigma ^k(N,\mathbb{R}^\ell)&=&{\rm dim}\ J^1(N,\mathbb{R}^\ell)- 
{\rm dim}\ \Sigma ^k(N,\mathbb{R}^\ell) \\
&=&(n-v+k)(\ell-v+k), 
\end{eqnarray*}
where $v={\rm min}\{n, \ell\}$. 
(\<For \> details on $\Sigma^k$ and $\Sigma^k(N,\mathbb{R}^\ell)$, see for example 
\cite{GG}, \<pp.\,\>60--61).
 
As applications of Theorem \ref{main}, we have 
the following Proposition \ref{submain}, 
Corollaries  
\ref{Morse function}, \ref{WU}, \ref{immersion} and \ref{corank}. 

\begin{proposition}\label{submain}
Let $N$ be a manifold of dimension $n$. 
Let $f$ be an immersion 
of $N$ into an open subset $U$ of $\mathbb{R}^m$. 
Let $F:U\to \mathbb{R}^\ell$ be a mapping. 
Then, there exists a subset $\Sigma$ of $\mathcal{L}(\mathbb{R}^{m},\mathbb{R}^{\ell})$ with Lebesgue measure zero  
such that for any $\pi \in \mathcal{L}(\mathbb{R}^{m},\mathbb{R}^{\ell})-\Sigma$, 
the mapping $j^1(F_\pi \circ f):
N\to J^1(N,\mathbb{R}^{\ell})$ is transverse to the submanifold 
$\Sigma ^k(N,\mathbb{R}^\ell)$ for any positive integer $k$ satisfying $1\leq k\leq v$. 
Especially, in the case of $\ell \geq2$, we have 
$k_0+1\leq v$ and it follows that the mapping $j^1(F_\pi \circ f)$ satisfies that 
$j^1(F_\pi \circ f)(N)\bigcap \Sigma^k(N,\mathbb{R}^\ell)=\emptyset$ 
for any positive integer $k$ satisfying $k_0+1\leq k \leq v$, 
where $k_0$ is the maximum integer 
satisfying 
$(n-v+k_0)(\ell-v+k_0)\leq n$ $(v={\rm min}\{n,\ell\})$. 
\end{proposition}

{\it Proof.}\qquad 
By Theorem \ref{main}, for any positive integer $k$ satisfying $1\leq k\leq v$, 
there exists a subset $\widetilde{\Sigma}_k$ of $\mathcal{L}(\mathbb{R}^{m},\mathbb{R}^{\ell})$ 
with Lebesgue measure zero 
such that for any $\pi \in \mathcal{L}(\mathbb{R}^{m},\mathbb{R}^{\ell})-\widetilde{\Sigma}_k$, 
the mapping $j^1(F_\pi \circ f):
N\to J^1(N,\mathbb{R}^{\ell})$ is transverse to $\Sigma^k(N,\mathbb{R}^\ell)$. 
Set $\Sigma=\bigcup_{k=1}^v\widetilde{\Sigma}_k$. 
Then, it is clearly seen that $\Sigma$ is a subset of 
$\mathcal{L}(\mathbb{R}^{m},\mathbb{R}^{\ell})$ with Lebesgue measure zero. 
Hence, it follows that there exists a subset $\Sigma$ of $\mathcal{L}(\mathbb{R}^{m},\mathbb{R}^{\ell})$ 
with Lebesgue measure zero  
such that for any $\pi \in \mathcal{L}(\mathbb{R}^{m},\mathbb{R}^{\ell})-\Sigma$, 
the mapping $j^1(F_\pi \circ f):
N\to J^1(N,\mathbb{R}^{\ell})$ is transverse to the submanifold 
$\Sigma ^k(N,\mathbb{R}^\ell)$ for any positive integer $k$ satisfying $1\leq k\leq v$. 

Now, we will consider the case of $\ell \geq 2$.  
Firstly, we will show that $k_0+1\leq v$ in the case. 
Suppose that $v\leq k_0$. 
Then, by $(n-v+k_0)(\ell-v+k_0)\leq n$, we have $n\ell \leq n$. 
This contradicts the assumption $\ell \geq 2$. 

Secondly, we will show that in the case of $\ell \geq 2$, 
the mapping $j^1(F_\pi \circ f): 
N\to J^1(N,\mathbb{R}^{\ell})$ satisfies that 
$j^1(F_\pi \circ f)(N)\bigcap \Sigma^k(N,\mathbb{R}^\ell)=\emptyset$ 
for any positive integer $k$ satisfying $k_0+1\leq k\leq v$. 
Suppose that there exist a positive integer $k$ $(k_0+1\leq k \leq v)$ 
and a point $q\in N$ such that 
$j^1(F_\pi \circ f)(q)\in \Sigma^k(N,\mathbb{R}^\ell)$. 
Since the mapping $j^1(F_\pi \circ f):
N\to J^1(N,\mathbb{R}^{\ell})$ is transverse to $\Sigma^k(N,\mathbb{R}^\ell)$ 
at the point $q$, the following holds: 
\begin{eqnarray*}
d(j^1(F_{\pi}\circ f))_{q}(T_{q}N)
+T_{j^1(F_{\pi}\circ f)(q)}\Sigma ^{k}(N,\mathbb{R}^{\ell})
=T_{j^1(F_{\pi}\circ f)(q)}J^1(N,\mathbb{R}^{\ell}).
\end{eqnarray*}
Hence, we have
\begin{eqnarray*}
{}&&{}{\rm dim}\ d(j^1(F_{\pi}\circ f))_{q}(T_{q}N)\\
&\geq &{\rm dim}\ T_{j^1(F_{\pi}\circ f)(q)}J^1(N,\mathbb{R}^{\ell})-
{\rm dim}\ T_{j^1(F_{\pi}\circ f)(q)}\Sigma ^{k}(N,\mathbb{R}^{\ell})\\
&=&{\rm codim}\ T_{j^1(F_{\pi}\circ f)(q)}\Sigma ^{k}(N,\mathbb{R}^{\ell}).
\end{eqnarray*}
Thus, we get $n\geq (n-v+k)(\ell-v+k)$. 
Since the given integer $k_0$ is the maximum integer 
satisfying 
$n\geq(n-v+k_0)(\ell-v+k_0)$, 
it follows that $k\leq k_0$. 
This contradicts the assumption $k_0+1\leq k$. 
\hfill\qed
\begin{remark}\label{remark k_0}
{\rm
\begin{enumerate}
\item
In Proposition \ref{submain}, by $(n-v+k_0)(\ell-v+k_0)\leq n$, it is clearly seen that  $k_0\geq0$.  
\item 
In Proposition \ref{submain}, in the case of $\ell=1$, we have $k_0+1>v$. 
Indeed, in the case, by $v=1$, 
we get $(n-1+k_0)k_0\leq n$. Hence, we have $k_0=1$.  
\end{enumerate}
}
\end{remark}

A mapping $g:N\to \mathbb{R}$ is called a {\it Morse function} 
if all of the singularities of the mapping $g$ are nondegenerate 
(for details on Morse function\<s\>, see for example, \cite{GG}, \<p.\,\>63). 
In the case of $(n,\ell)=(n,1)$, we have the following. 

\begin{corollary}\label{Morse function}
Let $N$ be a manifold of dimension $n$. 
Let $f$ be an immersion 
of $N$ into an open subset $U$ of $\mathbb{R}^m$. 
Let $F:U\to \mathbb{R}$ be a mapping. 
Then, there exists a subset $\Sigma$ of $\mathcal{L}(\mathbb{R}^{m},\mathbb{R})$  with Lebesgue measure zero 
such that for any $\pi \in \mathcal{L}(\mathbb{R}^{m},\mathbb{R})-\Sigma$, 
the mapping $F_\pi \circ f:
N\to \mathbb{R}$ is a Morse function. 
\end{corollary}
{\it Proof.}\qquad 
By Proposition \ref{submain}, 
there exists a subset $\Sigma$ with Lebesgue measure zero 
of $\mathcal{L}(\mathbb{R}^{m},\mathbb{R})$ 
such that for any $\pi \in \mathcal{L}(\mathbb{R}^{m},\mathbb{R})-\Sigma$, 
the mapping $j^1(F_\pi \circ f):
N\to J^1(N, \mathbb{R})$ is transverse to the submanifold $\Sigma ^1(N,\mathbb{R})$. 
Hence, if $q\in N$ is a singular point of the mapping $F_\pi \circ f$, then 
the point $q$ is nondegenerate. 
\hfill\qed
\par 
\bigskip 

For a given mapping $g:N\to \mathbb{R}^{2n-1}$ $(n\geq2)$, 
a singular point $q\in N$ is called a {\it singular point of Whitney umbrella} 
if there exist two germs of diffeomorphisms  
$H:(\mathbb{R}^{2n-1}, g(q)) \to (\mathbb{R}^{2n-1} , 0)$ and 
 $h:(N, q) \to (\R^n, 0)$ 
such that 
$H\circ g \circ h^{-1}(x_1,x_2,\ldots ,x_n)=(x_1^2,x_1x_2,\ldots ,x_1x_n, x_2, \ldots ,x_{n})$, 
where $(x_1,x_2,\ldots ,x_n)$ is a local coordinate around the point 
$h(q)=0\in \mathbb{R}^n$. 
In the case of $(n,\ell)=(n,2n-1)$ $(n\geq2)$, we have the following. 
\begin{corollary}\label{WU}
Let $N$ be a manifold of dimension $n$ $(n\geq2)$. 
Let $f$ be an immersion 
of $N$ into an open subset $U$ of $\mathbb{R}^m$. 
Let $F:U\to \mathbb{R}^{2n-1}$ be a mapping. 
Then, there exists a subset $\Sigma$ with Lebesgue measure zero 
of $\mathcal{L}(\mathbb{R}^{m},\mathbb{R}^{2n-1})$ 
such that for any $\pi \in \mathcal{L}(\mathbb{R}^{m},\mathbb{R}^{2n-1})-\Sigma$, 
any singular point of the mapping $F_\pi \circ f:
N\to \mathbb{R}^{2n-1}$ is 
a singular point of Whitney umbrella. 
\end{corollary}
{\it Proof.}\qquad 
By\<, \> for example, \cite{GG}, \<p.\,\>179, we see that 
a point $q\in N$ is a singular point of Whitney umbrella of the mapping $F_\pi \circ f$ 
if 
$j^1(F_\pi \circ f)(q)\in \Sigma^1(N,\mathbb{R}^{2n-1})$ and the mapping 
$j^1(F_\pi \circ f)$ is transverse to the submanifold 
$\Sigma^1(N,\mathbb{R}^{2n-1})$ at $q$. 
Set $\ell=2n-1$ and $v=n$ in Proposition \ref{submain}. 
Then, it is clearly seen that we have $k_0=1$ in Proposition \ref{submain}. 
Hence, there exists a subset $\Sigma$ of $\mathcal{L}(\mathbb{R}^{m},\mathbb{R}^{2n-1})$ 
with Lebesgue measure zero 
such that for any $\pi \in \mathcal{L}(\mathbb{R}^{m},\mathbb{R}^{2n-1})-\Sigma$, 
the mapping $F_\pi \circ f:
N\to \mathbb{R}^{2n-1}$ is transverse to 
$\Sigma ^k(N,\mathbb{R}^{2n-1})$ for any positive integer $k$ satisfying $1\leq k\leq n$, 
and the mapping satisfies that 
$j^1(F_\pi \circ f)(N)\bigcap \Sigma^k(N,\mathbb{R}^{2n-1})=\emptyset$ 
for any positive integer $k$ satisfying $2\leq k \leq n$. 
Thus, if a point $q\in N$ is a singular point of the mapping $F_\pi \circ f$, then 
it follows that 
$j^1(F_\pi \circ f)(q)\in \Sigma^1(N,\mathbb{R}^{2n-1})$ and
$j^1(F_\pi \circ f)$ is transverse to 
$\Sigma^1(N,\mathbb{R}^{2n-1})$ at $q$. 
\hfill\qed
\par 
\bigskip  

In the case of $\ell \geq2n$,  
the \<immersion \> property of a given mapping $f:N\to U$ 
is preserved by composing generic linearly perturbed mappings as follows:  
\begin{corollary}\label{immersion}
Let $N$ be a manifold of dimension $n$. 
Let $f$ be an immersion 
of $N$ into an open subset $U$ of $\mathbb{R}^{m}$. 
Let $F:U\to \mathbb{R}^{\ell}$ be a mapping $(\ell \geq2n)$. 
Then, there exists a subset $\Sigma$ of 
$\mathcal{L}(\mathbb{R}^{m},\mathbb{R}^{\ell})$ with Lebesgue measure zero 
such that for any $\pi \in \mathcal{L}(\mathbb{R}^{m},\mathbb{R}^{\ell})-\Sigma$, 
the mapping $F_\pi \circ f:
N\to \mathbb{R}^{\ell}$ is an immersion. 
\end{corollary}
{\it Proof.}\qquad 
It is clearly seen that the mapping $F_\pi \circ f:N\to \mathbb{R}^\ell$ 
is an immersion 
if and only if   
$j^1(F_\pi \circ f)(N)\bigcap \bigcup_{k=1}^n\Sigma^k(N,\mathbb{R}^\ell)=
 \emptyset$.
Set $v=n$ and $\ell \geq2n$ in Proposition \ref{submain}. 
Then, it is clearly seen that $k_0\leq0$. By Remark \ref{remark k_0}, 
we get $k_0=0$. 
Hence, 
there exists a subset $\Sigma$ of $\mathcal{L}(\mathbb{R}^{m},\mathbb{R}^{\ell})$ with Lebesgue measure zero 
such that for any $\pi \in \mathcal{L}(\mathbb{R}^{m},\mathbb{R}^{\ell})-\Sigma$, 
the mapping $j^1(F_\pi \circ f):
N\to J^1(N,\mathbb{R}^{\ell})$ satisfies that 
$j^1(F_\pi \circ f)(N)\bigcap \Sigma^k(N,\mathbb{R}^\ell)=\emptyset$ 
for any positive integer $k$ $(1\leq k \leq n)$. 
\hfill\qed
\par 
\bigskip 
A mapping $g:N\to \mathbb{R}^\ell$ 
{\it has corank  at most k singular points} 
if  
\begin{eqnarray*}
\sup\left\{{\rm corank\ }dg_q\mid q\in N \right\}\leq k, 
\end{eqnarray*}
where ${\rm corank\ }dg_q={\rm min} \{n,\ell \}-{\rm rank\ } dg_q$. 
By Proposition \ref{submain}, we have the following corollary.  
\begin{corollary}\label{corank}
Let $N$ be a manifold of dimension $n$. 
Let $f$ be an immersion 
of $N$ into an open subset $U$ of $\mathbb{R}^m$. 
Let $F:U\to \mathbb{R}^\ell$ be a mapping. 
Let $k_0$ be the maximum integer 
satisfying 
$(n-v+k_0)(\ell-v+k_0)\leq n$ $(v={\rm min}\{n,\ell\})$.
Then, there exists a subset $\Sigma$ of $\mathcal{L}(\mathbb{R}^{m},\mathbb{R}^{\ell})$ 
with Lebesgue measure zero  
such that for any $\pi \in \mathcal{L}(\mathbb{R}^{m},\mathbb{R}^{\ell})-\Sigma$, 
the mapping $F_\pi \circ f:
N\to \mathbb{R}^\ell$ has corank at most $k_0$ singular points. 
\end{corollary}
\par 
\bigskip 

\subsection{Applications of Theorem \ref{main2}}\label{application2}
\begin{proposition}\label{normal}
Let $N$ be a manifold of dimension $n$. 
Let $f$ be an injection of $N$ into an open subset  
$U$ of $\mathbb{R}^m$. 
Let $F:U\to \mathbb{R}^\ell$ be a mapping. 
If $(s_f-1)\ell>n s_f$, then 
there exists a subset $\Sigma$ of $\mathcal{L}(\mathbb{R}^m, \mathbb{R}^\ell)$ 
with Lebesgue measure zero such that 
for any $\pi \in \mathcal{L}(\mathbb{R}^m, \mathbb{R}^\ell)-\Sigma $, 
$F_\pi \circ f:N\to \mathbb{R}^\ell$ is a mapping with normal crossings 
satisfying $(F_\pi \circ f)^{(s_f)}(N^{(s_f)})\bigcap\Delta _{s_f}=\emptyset$. 
\end{proposition}
{\it Proof.}\qquad 
By Theorem \ref{main2}, 
there exists a subset $\Sigma$ of $\mathcal{L}(\mathbb{R}^m, \mathbb{R}^\ell)$ with Lebesgue measure zero such that 
for any $\pi \in \mathcal{L}(\mathbb{R}^m, \mathbb{R}^\ell)-\Sigma $, 
\<and \> 
for any $s$ $(2\leq s \leq s_f)$, 
the mapping $(F_\pi \circ f)^{(s)}:N^{(s)}\to (\R^\ell)^s$ is transverse to the 
submanifold $\Delta_s$. 
Hence, in order to show Proposition \ref{normal},  
it is sufficient to show that 
for any $\pi \in \mathcal{L}(\mathbb{R}^m, \mathbb{R}^\ell)-\Sigma $, 
the mapping $(F_\pi \circ f)^{(s_f)}$ satisfies that 
$(F_\pi \circ f)^{(s_f)}(N^{(s_f)})\bigcap \Delta_{s_f}=\emptyset$. 

Suppose that there exists an element 
$\pi \in \mathcal{L}(\mathbb{R}^m, \mathbb{R}^\ell)-\Sigma $ 
such that there exists a point $q\in N^{(s_f)}$ 
satisfying $(F_{\pi}\circ f)^{(s_f)}(q)\in \Delta_{s_f}$. 
Since $(F_{\pi}\circ f)^{(s_f)}$ is transverse to $\Delta_{s_f}$, 
we have the following: 
\begin{eqnarray*}
d((F_{\pi}\circ f)^{(s_f)})_{q}(T_{q}N^{(s_f)})
+T_{(F_{\pi}\circ f)^{(s_f)}(q)}\Delta_{s_f}
=T_{(F_{\pi}\circ f)^{(s_f)}(q)}(\mathbb{R}^\ell)^{s_f}.
\end{eqnarray*}
Hence, we have
\begin{eqnarray*}
{}&&{}{\rm dim}\ d((F_{\pi}\circ f)^{(s_f)})_{q}(T_{q}N^{(s_f)})\\
&\geq &{\rm dim}\ T_{(F_{\pi}\circ f)^{(s_f)}(q)}(\mathbb{R}^\ell)^{s_f}-
{\rm dim}\ T_{(F_{\pi}\circ f)^{(s_f)}(q)}\Delta_{s_f}\\
&=&{\rm codim}\ T_{(F_{\pi}\circ f)^{(s_{f})}(q)}\Delta_{s_f}.
\end{eqnarray*}
Thus, we get $ns_f\geq (s_f-1)\ell$. 
This contradicts the assumption $(s_f-1)\ell>ns_f$. 
\hfill\qed
\par 
\bigskip  
In the case of $\ell >2n$, 
the \<injection \> property of a given mapping $f:N\to U$ 
is preserved by composing generic linearly perturbed mappings as follows:  
\begin{corollary}\label{injective}
Let $N$ be a manifold of dimension $n$. 
Let $f$ be an injection of $N$ into an open subset  
$U$ of $\mathbb{R}^m$. 
Let $F:U\to \mathbb{R}^\ell$ be a mapping. 
If $\ell>2n$, then 
there exists a subset $\Sigma$ of $\mathcal{L}(\mathbb{R}^m, \mathbb{R}^\ell)$ 
with Lebesgue measure zero such that 
for any $\pi \in \mathcal{L}(\mathbb{R}^m, \mathbb{R}^\ell)-\Sigma $, 
the mapping $F_\pi \circ f:N\to \mathbb{R}^\ell$ is injective. 
\end{corollary}
{\it Proof.}\qquad 
Since $s_f\geq 2$ and $\ell>2n$, it is easily seen that the dimension pair $(n, \ell)$ 
satisfies the assumption $(s_f-1)\ell>ns_f$ of Proposition \ref{normal}. 
Indeed, \<from \> $\ell>2n$, it follows that $(s_f-1)\ell>2n(s_f-1)$. 
By $s_f\geq 2$, we get $2n(s_f-1)\geq ns_f$. 
 
Hence, 
by Proposition \ref{normal}, 
there exists a subset $\Sigma$ of $\mathcal{L}(\mathbb{R}^m, \mathbb{R}^\ell)$ 
with Lebesgue measure zero such that 
for any $\pi \in \mathcal{L}(\mathbb{R}^m, \mathbb{R}^\ell)-\Sigma $, 
the mapping $(F_\pi \circ f)^{(2)} : N^{(2)} \to (\mathbb{R}^\ell)^2$ is 
transverse to $\Delta_2$. 
In order to show Corollary \ref{injective}, 
it is sufficient to show that 
the mapping $(F_\pi \circ f)^{(2)}$ satisfies that 
$(F_\pi \circ f)^{(2)}(N^{(2)})\bigcap \Delta_2 =\emptyset$. 

Suppose that there exists a point $q\in N^{(2)}$ such that 
$(F_\pi \circ f)^{(2)}(q)\in \Delta_2$. 
Then, we have the following: 
\begin{eqnarray*}
d((F_{\pi}\circ f)^{(2)})_{q}(T_{q}N^{(2)})
+T_{(F_{\pi}\circ f)^{(2)}(q)}\Delta_{2}
=T_{(F_{\pi}\circ f)^{(2)}(q)}(\mathbb{R}^\ell)^{2}.
\end{eqnarray*}
Hence, we have
\begin{eqnarray*}
{}&&{}{\rm dim}\ d((F_{\pi}\circ f)^{(2)})_{q}(T_{q}N^{(2)})\\
&\geq &{\rm dim}\ T_{(F_{\pi}\circ f)^{(2)}(q)}(\mathbb{R}^\ell)^{2}-
{\rm dim}\ T_{(F_{\pi}\circ f)^{(2)}(q)}\Delta_{2}\\
&=&{\rm codim}\ T_{(F_{\pi}\circ f)^{(2)}(q)}\Delta_{2}.
\end{eqnarray*}
Thus, we get $2n\geq \ell$. 
This contradicts the assumption $\ell>2n$. 
\hfill\qed
\par 
\bigskip  

By combining Corollaries \ref{immersion} and \ref{injective}, 
we have the following. 
\begin{corollary}\label{injective immersion}
Let $N$ be a manifold of dimension $n$. 
Let $f$ be an injective immersion of $N$ into an open subset  
$U$ of $\mathbb{R}^m$. 
Let $F:U\to \mathbb{R}^\ell$ be a mapping. 
If $\ell>2n$, then 
there exists a subset $\Sigma$ of $\mathcal{L}(\mathbb{R}^m, \mathbb{R}^\ell)$ 
with Lebesgue measure zero such that 
for any $\pi \in \mathcal{L}(\mathbb{R}^m, \mathbb{R}^\ell)-\Sigma $, 
the mapping $F_\pi \circ f:N\to \mathbb{R}^\ell$ is an injective immersion. 
\end{corollary}

In Corollary \ref{injective immersion}, suppose that the mapping 
$F_\pi \circ f:N\to \mathbb{R}^\ell$ is proper. 
Then, an injective immersion  $F_\pi \circ f$ 
is necessarily an embedding 
(see \cite{GG}, p.\,11). Thus, we get the following. 
\begin{corollary}\label{embedding}
Let $N$ be a compact manifold of dimension $n$. 
Let $f$ be an embedding of $N$ into an open subset  
$U$ of $\mathbb{R}^m$. 
Let $F:U\to \mathbb{R}^\ell$ be a mapping. 
If $\ell>2n$, then 
there exists a subset $\Sigma$ of $\mathcal{L}(\mathbb{R}^m, \mathbb{R}^\ell)$ 
with Lebesgue measure zero such that 
for any $\pi \in \mathcal{L}(\mathbb{R}^m, \mathbb{R}^\ell)-\Sigma $, 
the mapping $F_\pi \circ f:N\to \mathbb{R}^\ell$ is an embedding. 
\end{corollary}
\section{Further applications}\label{section 6}
\subsection{Introduction of generalized distance-squared mappings}

Let $p_i=(p_{i1}, p_{i2}, \ldots, p_{im})$  $(1\le i\le \ell)$ 
(resp., $A=(a_{ij})_{1\le i\le \ell, 1\le j\le m}$) 
be points of $\mathbb{R}^m$ 
(resp., an $\ell\times m$ matrix with all entries being non-zero real numbers). 
Set $p=(p_1,p_2,\ldots,p_{\ell})\in (\mathbb{R}^m)^{\ell}$. 
Let $G_{(p, A)}:\mathbb{R}^m \to \mathbb{R}^\ell$ be the mapping 
defined by 
{\small 
\[
G_{(p, A)}(x)=\left(
\sum_{j=1}^m a_{1j}(x_j-p_{1j})^2, 
\sum_{j=1}^m a_{2j}(x_j-p_{2j})^2, 
\ldots, 
\sum_{j=1}^m a_{\ell j}(x_j-p_{\ell j})^2
\right), 
\] }where $x=(x_1, x_2, \ldots, x_m)\in \mathbb{R}^m$. 
The mapping $G_{(p, A)}$ is called a {\it generalized distance-squared mapping}, 
and the $\ell$-tuple of points 
$p=(p_1,p_2,\ldots ,p_{\ell}) \in (\mathbb{R}^m)^{\ell}$ 
is called the {\it central point} 
of the generalized distance-squared mapping $G_{(p,A)}$. 
A {\it distance-squared mapping} $D_p$ 
(resp., {\it Lorentzian distance-squared mapping}
$L_p$) is the mapping $G_{(p,A)}$ 
satisfying that each entry of $A$ is equal to $1$ 
(resp., $a_{i1}=-1$ and $a_{ij}=1$ $(j\ne 1)$). 

In \cite{D} (resp., \cite{L}), 
a classification result of distance-squared mappings 
(resp., Lorentzian distance-squared mappings) is given. 

In \cite{G1}, 
a classification result of 
generalized distance-squared mappings of the plane into the plane 
is given. 
If the rank of $A$ is equal to two, 
then 
a generalized distance-squared mapping 
having a generic central point 
is a mapping of which any singular point is 
a {\rm fold point} 
except one {\rm cusp point}.     
The singular set is a rectangular hyperbola. 
If the rank of $A$ is equal to one, then 
a generalized distance-squared mapping 
having a generic central point is $\mathcal{A}$-equivalent 
to {\rm the normal form of fold singularity} $(x_1,x_2)
\mapsto (x_1,x_2^2)$. 

In \cite{G2}, 
a classification result of 
generalized distance-squared mappings 
of $\mathbb{R}^{m+1}$ into $\mathbb{R}^{2m+1}$ 
is given. 
If the rank of $A$ is equal to $m+1$, 
then  
a generalized distance-squared mapping 
having a generic central point 
is $\mathcal{A}$-equivalent to 
{\rm the normal form of Whitney umbrella} 
$
(x_1,\<x_2,\>\ldots ,x_{m+1})\mapsto 
(x_1^2,x_1x_2,\ldots ,x_1x_{m+1},x_2,\ldots ,x_{m+1})
$. 
If the rank of $A$ is strictly smaller than $m+1$, then 
a generalized distance-squared mapping 
having a generic central point 
is $\mathcal{A}$-equivalent to the inclusion 
$
(x_1,x_2,\ldots ,x_{m+1})\mapsto (x_1,x_2,\ldots ,x_{m+1},0,\ldots ,0)
$. 
\par 
Namely, in \cite{D}, \cite{L}, \cite{G2} and \cite{G1}, 
the properties of generic generalized distance-squared mappings 
are investigated. 
Hence, it is natural to investigate 
the properties of compositions with generic generalized distance-squared mappings. 

We have another original motivation. 
Height functions and distance-squared functions have been investigated 
in detail so far, 
and they are useful tools   
in the applications of singularity theory to differential geometry 
(for instance, see \cite{CS}). 
A mapping in which each component is a height function is 
nothing but a projection. 
Projections as well as height functions or distance-squared functions 
have been investigated so far.  
In \cite{GP}, compositions of 
generic projections and embeddings are investigated.

On the other hand, 
a mapping in which each component 
is a distance-squared function is a distance-squared mapping. 
In addition, the notion of a generalized distance-squared mapping is 
an extension of 
that of a distance-squared mapping. 
Therefore, it is natural to investigate compositions 
with generic generalized distance-squared mappings as well as projections. 
\par 

\subsection{Applications of Theorem \ref{main} to $G_{(p, A)}:\mathbb{R}^m \to \mathbb{R}^\ell$}
\begin{proposition} 
\label{Gmain}
Let $N$ be a manifold of dimension $n$. 
Let $f:N\to \mathbb{R}^m$ be an immersion. 
Let $A=(a_{ij})_{1\leq i \leq \ell, 1\leq j \leq m}$ 
be an $\ell \times m$ matrix with all entries being non-zero real numbers. 
If $X$ is an $\mathcal{A}^1$-invariant submanifold of $J^1(n,\ell)$, 
then there exists a subset $\Sigma$ of $(\mathbb{R}^m)^\ell$ 
with Lebesgue measure zero  
such that for any $p \in (\mathbb{R}^m)^\ell-\Sigma$, 
the mapping $j^1(G_{(p,A)} \circ f):
N\to J^1(N,\mathbb{R}^\ell)$ is transverse to 
the submanifold $X(N,\mathbb{R}^\ell)$.
\end{proposition}
{\it Proof.}\qquad 
Let $H:\mathbb{R}^{\ell} \to \mathbb{R}^{\ell}$ 
be a diffeomorphism of the target for deleting constant terms. 
The composition $H\circ G_{(p,A)}:\mathbb{R}^{m}\to \mathbb{R}^{\ell}$ 
is given as follows:
\begin{eqnarray*}
H\circ G_{(p, A)}(x)&=&\left(
\sum_{j=1}^m a_{1j}x_j^2-2\sum_{j=1}^m a_{1j}p_{1j}x_j, 
\sum_{j=1}^m a_{2j}x_j^2-2\sum_{j=1}^m a_{2j}p_{2j}x_j, \right. 
\\&&
\left. 
\ldots, 
\sum_{j=1}^m a_{\ell j}x_j^2-2\sum_{j=1}^m a_{\ell j}p_{\ell j}x_j
\right), 
\end{eqnarray*}
where $x=(x_1,x_2,\ldots ,x_m)$. 

Let $\psi :(\mathbb{R}^m)^\ell \to \mathcal{L}(\mathbb{R}^{m},\mathbb{R}^{\ell})$ 
be the mapping defined by 
\begin{eqnarray*}
\psi (p_{11},p_{12},\ldots ,p_{\ell m})=-2(a_{11}p_{11},a_{12}p_{12},\ldots ,
a_{\ell m}p_{\ell m}). 
\end{eqnarray*}
Remark that we have the natural identification 
$\mathcal{L}(\mathbb{R}^{m},\mathbb{R}^{\ell})=(\mathbb{R}^m)^\ell$. 
Since $a_{ij}\not =0$ for any $i$, $j$ $(1\leq i \leq \ell$, $1\leq j \leq m)$, 
it is clearly seen that $\psi$ is a $C^\infty$ diffeomorphism. 

Set $F_{i}(x)=\sum_{j=1}^m a_{ij}x_j^2$ $(1\leq i \leq \ell)$ and 
$F=(F_1,F_2,\ldots ,F_\ell)$.  
By Theorem \ref{main}, 
there exists a subset $\Sigma$ of $\mathcal{L}(\mathbb{R}^{m},\mathbb{R}^{\ell})$ with Lebesgue measure zero  
such that 
for any $\pi \in \mathcal{L}(\mathbb{R}^{m},\mathbb{R}^{\ell})-\Sigma$, 
the mapping $j^{1}(F_{\pi}\circ f):N\to J^{1}(N,\mathbb{R}^\ell)$ 
is transverse to $X(N,\mathbb{R}^\ell)$. 
Since $\psi^{-1}:\mathcal{L}(\mathbb{R}^{m},\mathbb{R}^{\ell}) \to 
(\mathbb{R}^m)^\ell$ is a $C^\infty$ mapping, 
$\psi^{-1}(\Sigma)$ is a subset of $(\mathbb{R}^m)^\ell$ 
with Lebesgue measure zero. 
For any $p\in (\mathbb{R}^m)^\ell-\psi^{-1}(\Sigma)$, 
we have $\psi(p)\in \mathcal{L}(\mathbb{R}^{m},\mathbb{R}^{\ell})-\Sigma$. 
Hence, for any $p\in (\mathbb{R}^m)^\ell-\psi^{-1}(\Sigma)$, 
the mapping 
$j^{1}(H\circ G_{(p, A)}\circ f):N\to J^{1}(N,\mathbb{R}^\ell)$ 
is transverse to $X(N,\mathbb{R}^\ell)$. 
Then, since $H:\mathbb{R}^{\ell} \to \mathbb{R}^{\ell}$ 
is a diffeomorphism, 
the mapping 
$j^{1}(G_{(p, A)}\circ f):N\to J^{1}(N,\mathbb{R}^\ell)$ 
is transverse to $X(N,\mathbb{R}^\ell)$. 
\hfill\qed

\begin{remark}\label{remarkGmain}
{\rm 
As applications of Proposition \ref{Gmain}, regarding generalized distance-squared mappings, 
we get analogies of  Proposition \ref{submain}, 
Corollaries \ref{Morse function}, \ref{WU}, \ref{immersion} and \ref{corank}.  
}
\end{remark}
\subsection{Applications of Theorem \ref{main2} to $G_{(p, A)}:\mathbb{R}^m \to \mathbb{R}^\ell$}
By Theorem~\ref{main2}, 
we get the following proposition, 
which can be proved by the same argument as in the proof 
of Proposition~\ref{Gmain}, and we omit the proof. 
\begin{proposition}\label{Gmain2}
Let $N$ be a manifold of dimension $n$. 
Let $f : N\to \mathbb{R}^m$ be an injection. 
Let $A=(a_{ij})_{1\leq i \leq \ell, 1\leq j \leq m}$ 
be an $\ell \times m$ matrix with all entries being non-zero real numbers. 
Then, there exists a subset $\Sigma$ of $(\mathbb{R}^m)^\ell$ with Lebesgue measure zero such that 
for any $p \in (\mathbb{R}^m)^\ell-\Sigma $, 
and 
for any $s$ $(2\leq s \leq s_f)$, 
the mapping $(G_{(p,A)} \circ f)^{(s)}:N^{(s)}\to (\R^\ell)^s$ is transverse to the 
submanifold $\Delta_s$. 
Moreover, if the mapping $G_{(p,A)}$ satisfies that $| G_{(p,A)}^{-1}(y) | \leq s_f$ for any $y\in\mathbb{R}^\ell$, 
then $G_{(p,A)} \circ f:N\to \mathbb{R}^\ell$ is a mapping with normal crossings.    
\end{proposition}
\begin{remark}\label{remarkGmain2}
{\rm 
As applications of Proposition~\ref{Gmain2}, regarding generalized distance-squared mappings, 
we get analogies of Proposition \ref{normal}, 
Corollaries \ref{injective}, \ref{injective immersion} and \ref{embedding}. 
}
\end{remark}

As the special case of the classification result of 
distance squared mappings (resp., Lorentzian distance-squared mappings) in \cite{D} (resp., \cite{L}), 
we have Lemma \ref{Dlemma}. 
\begin{lemma}[\cite{D}, \cite{L}] \label{Dlemma}
We have the following. 
\begin{enumerate}
\item
For any $p \in \mathbb{R}$, 
the mappings $D_p: \mathbb{R}\to \mathbb{R}$ and 
$L_p: \mathbb{R}\to \mathbb{R}$
are $\mathcal{A}$-equivalent to $x\mapsto x^2$. 
\item
For $m \geq 2$, 
there exists a subset $\Sigma_D$ $($resp., $\Sigma_L$$)$ of 
$(\mathbb{R}^m)^m$ with Lebesgue measure zero such that 
for any $p \in (\mathbb{R}^m)^m-\Sigma_D$ 
$($resp., $p \in (\mathbb{R}^m)^m-\Sigma_L$$)$, 
the mapping $D_p: \mathbb{R}^m\to \mathbb{R}^m$ 
$($resp., $L_p: \mathbb{R}^m\to \mathbb{R}^m$$)$
is $\mathcal{A}$-equivalent to the normal form of definite fold mappings 
$(x_1, x_2, \ldots ,x_m)\mapsto (x_1, x_2, \ldots ,x_{m-1},x_m^2)$. 
\item
In the case of $1\leq m< \ell$, 
there exists a subset $\Sigma_D$ $($resp., $\Sigma_L$$)$ 
of $(\mathbb{R}^m)^\ell$ with Lebesgue measure zero such that 
for any $p \in (\mathbb{R}^m)^\ell-\Sigma_D$ 
$($resp., $p \in (\mathbb{R}^m)^\ell-\Sigma_L$$)$, 
the mapping $D_p: \mathbb{R}^m\to \mathbb{R}^\ell$ 
$($resp., $L_p: \mathbb{R}^m\to \mathbb{R}^\ell$$)$ 
is $\mathcal{A}$-equivalent to the inclusion 
$(x_1,x_2, \ldots ,x_m)\mapsto (x_1,x_2, \ldots ,x_{m},0,\ldots ,0)$. 
\end{enumerate}
\end{lemma}

\begin{proposition}\label{DL}
Let $N$ be a manifold of dimension $n$. 
Let $f  : N\to \mathbb{R}^m$ be an injection. 
Then, the following holds:
\begin{enumerate}
\item
For $m \geq 1$, 
there exists a subset $\Sigma_D$ $($resp., $\Sigma_L$$)$ of 
$(\mathbb{R}^m)^m$ with Lebesgue measure zero such that 
for any $p \in (\mathbb{R}^m)^m-\Sigma_D$ 
$($resp., $p \in (\mathbb{R}^m)^m-\Sigma_L$$)$, 
$D_p \circ f:N\to \mathbb{R}^m$ 
$($resp., $L_p \circ f:N\to \mathbb{R}^m$$)$ 
is a mapping with normal crossings. 
\item
In the case of $1 \leq m<\ell$, 
there exists a subset $\Sigma_D$ $($resp., $\Sigma_L$$)$ of 
$(\mathbb{R}^m)^\ell$ with Lebesgue measure zero such that 
for any $p \in (\mathbb{R}^m)^\ell-\Sigma_D$ 
$($resp., $p \in (\mathbb{R}^m)^\ell-\Sigma_L$$)$, 
the mapping $D_p \circ f:N\to \mathbb{R}^\ell$ 
$($resp., $L_p \circ f:N\to \mathbb{R}^\ell$$)$ 
is an injection. 
\end{enumerate}
\end{proposition}
{\it Proof.}\qquad 
The proof for distance-squared mappings is the same 
as that for Lorentzian distance-squared mappings. 
Hence, it is sufficient to give the proof for distance-squared mappings. 

Firstly, we will show the assertion 1. 
From Lemma \ref{Dlemma}, 
there exists a subset $\Sigma_1$ of $(\mathbb{R}^m)^m$ with Lebesgue measure zero such that 
for any $p \in (\mathbb{R}^m)^m-\Sigma _1$, 
the mapping $D_p: \mathbb{R}^m\to \mathbb{R}^m$ 
satisfies that $|D_p^{-1}(y)| \leq 2$ 
for any $y\in \mathbb{R}^m$. On the other hand, from Proposition \ref{Gmain2}, 
there exists a subset $\Sigma_2$ of $(\mathbb{R}^m)^m$ with Lebesgue measure zero such that 
for any $p \in (\mathbb{R}^m)^m-\Sigma_2 $, 
if $D_p$ satisfies that $| D_p^{-1}(y) | \leq s_f$ for any $y\in\mathbb{R}^m$, 
then $D_p \circ f:N\to \mathbb{R}^m$ 
is a mapping with normal crossings. 
Set $\Sigma_D=\Sigma_1\cup \Sigma_2$. 
It is clearly seen that $\Sigma_D$ is a subset of 
 $(\mathbb{R}^m)^m$ with Lebesgue measure zero. 
Then, 
 for any $p \in (\mathbb{R}^m)^m-\Sigma_D$, $D_p \circ f:N\to \mathbb{R}^m$ 
is a mapping with normal crossings. 

In the case of $m<\ell$, since from Lemma \ref{Dlemma}, 
there exists a subset $\Sigma_D$ of 
$(\mathbb{R}^m)^\ell$ with Lebesgue measure zero such that 
for any $p \in (\mathbb{R}^m)^\ell-\Sigma_D$, 
the mapping $D_p : \mathbb{R}^m\to \R^\ell$ 
is $\mathcal{A}$-equivalent to the inclusion, the assertion 2 holds. 
\hfill\qed
 
\bigskip
By combining Proposition \ref{DL} and the analogy of 
Corollary \ref{immersion} in Remark \ref{remarkGmain}, we have the 
following. 
\begin{corollary}\label{immersion with}
Let $N$ be a manifold of dimension $n$. 
Let $f : N\to \mathbb{R}^m$ be an injective immersion $(2n\leq m)$. 
Then, 
there exists a subset $\Sigma_D$ $($resp., $\Sigma_L$$)$ 
of $(\mathbb{R}^{m})^{m}$ 
with Lebesgue measure zero such that 
for any $p \in (\mathbb{R}^{m})^{m}-\Sigma_D$ 
$($resp., $p \in (\mathbb{R}^{m})^{m}-\Sigma_L$$)$, 
the mapping $D_p \circ f:N\to \mathbb{R}^{m}$ 
$($resp., $L_p \circ f:N\to \mathbb{R}^{m}$$)$ 
is an immersion with normal crossings. 
\end{corollary}
In Corollary \ref{immersion with}, if $m=2n$ and the mapping 
$D_p \circ f:N\to \mathbb{R}^{2n}$ (resp., $L_p \circ f:N\to \mathbb{R}^{2n}$) 
is proper, 
then the immersion with normal crossings $D_p \circ f:N\to \mathbb{R}^{2n}$ 
(resp., $L_p \circ f:N\to \mathbb{R}^{2n}$) 
is necessarily stable 
(see \cite{GG}, p.\,86).   
Thus, we get the following. 
\begin{corollary}\label{stable}
Let $N$ be a compact manifold of dimension $n$. 
Let $f : N\to \mathbb{R}^{2n}$ be an embedding. 
Then, 
there exists a subset $\Sigma_D$ $($resp., $\Sigma_L$$)$ 
of $(\mathbb{R}^{2n})^{2n}$ with Lebesgue measure zero such that 
for any $p \in (\mathbb{R}^{2n})^{2n}-\Sigma_D$ 
$($resp., $p \in (\mathbb{R}^{2n})^{2n}-\Sigma_L$$)$, 
the mapping $D_p \circ f:N\to \mathbb{R}^{2n}$ 
$($resp., $L_p \circ f:N\to \mathbb{R}^{2n}$$)$
is stable. 
\end{corollary}
Remark that the dimension of the target space in Corollary \ref{stable} is 
smaller than that in Corollary \ref{embedding}. 

\section{Appendix}\label{section 7}
In this section, 
the main theorems in \cite{F1} and \cite{GP} are stated. 
For this, 
we prepare some notions. 

Let $N$ and $P$ be manifolds. 
Let ${}_s J^r(N,P)$ be the space consisting of 
elements $(j^r g(q_1),j^r g(q_2), \ldots ,j^r g(q_s))\in J^r(N,P)^s$ 
satisfying $(q_1,q_2,\ldots ,q_s)\in N^{(s)}$. 
Since $N^{(s)}$ is an open submanifold of $N^s$, 
the space ${}_s J^r(N,P)$ is also an open submanifold of $J^r(N,P)^s$.
For a given mapping $g:N\to P$, 
the mapping ${}_s j^rg:N^{(s)}\to {}_s J^r(N,P)$ 
is defined by 
\< 
$(q_1,q_2,\ldots ,q_s) \mapsto (j^rg(q_1),j^rg(q_2),\ldots ,j^rg(q_s))\>$.  

Let $W$ be a submanifold of ${}_s J^r(N,P)$. 
A mapping $g:N\to P$ will be said to be {\it transverse with respect to} $W$ 
if ${}_s j^r g:N^{(s)}\to {}_s J^r(N,P)$ is transverse to $W$. 

Following Mather (\cite{GP}), we can partition $P^s$ as follows. 
Given any partition $\Pi$ of $\{1,2, \ldots ,s\}$,  
let $P^{\Pi}$ denote the set of $s$-tuples 
$(y_1,y_2,\ldots ,y_s)\in P^s$ 
such that $y_i=y_j$ 
if and only if the 
two positive integers 
$i$ and $j$ are in the same member of the partition $\Pi$. 

Let Diff $N$ denote 
the group of diffeomorphisms of $N$. 
We have the natural action of Diff $N$ $\times $ Diff $P$ 
on ${}_s J^r(N,P)$ 
such that for a mapping $g:N\to P$, 
the equality 
\<$(h,H)\cdot {}_sj^rg(q)={}_sj^r(H\circ g\circ h^{-1})(q')$ 
holds, where $q=(q_1,q_2, \ldots ,q_s)$ and $q'=(h(q_1),h(q_2), \ldots ,h(q_s))\>$. 
A subset $W$ of ${}_s J^r(N,P)$ is said to be {\it invariant} 
if it is invariant under this action. 

We recall the following identification (7.1) from \cite{GP}. 
For $q=(q_1,q_2, \ldots ,q_s)\in N^{(s)}$, 
let $g:U\to P$ be a mapping defined 
in a neighborhood $U$ of \<$\{q_1,q_2, \ldots ,q_s\}\>$ in $N$, and 
let $z={}_s j^r g(q)$, \<$q'=(g(q_1),g(q_2), \ldots ,g(q_s))\>$. 
Let ${}_s J^r(N,P)_q$ and ${}_s J^r(N,P)_{q,q'}$ denote 
the fibers of ${}_s J^r(N,P)$ over $q$ and over $(q,q')$ respectively. 
Let $J^r(N)_q$ denote the $\mathbb{R}$-algebra of 
$r$-jets at $q$ of functions on $N$. 
Namely, 
\begin{eqnarray*}
J^r(N)_q={}_s J^r(N,\mathbb{R})_q.
\end{eqnarray*} 
Set $g^*TP=\bigcup_{\widetilde{q}\in U}T_{g(\widetilde{q})}P$, 
where $TP$ is the tangent bundle of $P$.    
Let $J^r(g^*TP)_q$ denote 
the $J^r(N)_q$-module of $r$-jets at $q$ of 
sections of the bundle $g^*TP$. 
Let $\mathfrak{m}_q$ be the ideal in $J^r(N)_q$ 
consisting of jets of functions which vanish at $q$. 
Namely,   
\begin{eqnarray*}
\mathfrak{m}_q=\{{}_s j^r h(q)\in {}_s J^r(N,\mathbb{R})_q \mid 
h(q_1)=h(q_2)=\cdots =h(q_s)=0\}. 
\end{eqnarray*}
Let $\mathfrak{m}_q J^r(g^*TP)_q$ be the set consisting of finite sums 
of products of 
an element of $\mathfrak{m}_q$ and an element of $J^r(g^*TP)_q$. 
Namely, we set  
\begin{eqnarray*}
\mathfrak{m}_q J^r(g^*TP)_q
=J^r(g^*TP)_q\cap \{{}_sj^r \xi(q)\in {}_s J^r(N,TP)_q \mid \xi(q_1)=
\xi(q_2)=
\cdots =\xi(q_s)=0\}.
\end{eqnarray*}
Then, it is easily seen that 
we have the following canonical identification of 
$\R$-vector spaces: 
\[
T({}_s J^r(N,P)_{q,q'})_z=\mathfrak{m}_q J^r(g^*TP)_q.  \eqno (7.1)
\]
Let $W$ be a non-empty submanifold of 
${}_s J^r(N,P)$. 
Choose $q=(q_1,q_2, \ldots ,q_s)\in N^{(s)}$ and 
$g:N\to P$, and 
set $z={}_s j^r g(q)$ and $q'=(g(q_1),g(q_2),\ldots ,g(q_s))$. 
Suppose that the choice is made so that $z\in W$. 
Set $W_{q,q'}=\widetilde{\pi}^{-1}(q,q')$, where 
$\widetilde{\pi}:W\to N^{(s)}\times P^s$ is defined by 
$\widetilde{\pi}({}_s j^r \widetilde{g}(\widetilde{q}))
=(\widetilde{q}, (
\widetilde{g}(\widetilde{q}_1), \widetilde{g}(\widetilde{q}_2), \ldots ,
\widetilde{g}(\widetilde{q}_s)))$ and 
$\widetilde{q}=(\widetilde{q}_1, \widetilde{q}_2, \ldots , \widetilde{q}_s)\in N^{(s)}$.  

Then, under the identification (7.1), the tangent space 
$T(W_{q,q'})_z$ can be identified with a vector subspace of 
$\mathfrak{m}_q$$J^r(g^*TP)_q$. 
We denote this vector subspace by $E(g,q,W)$. 

\begin{definition}\label{modular}
{\rm 
The submanifold $W$ is said to be {\it modular} 
if conditions $(\alpha)$ and $(\beta)$ below are satisfied.
\begin{enumerate}
\item[$(\alpha)$]
The set $W$ is an invariant submanifold of ${}_s J^r(N,P)$, 
and lies over $P^{\Pi}$ 
for some partition $\Pi$ of $\{1,2, \ldots ,s\}$.
\item[$(\beta)$] 
For any $q\in N^{(s)}$ and any mapping $g:N \to P$ 
such that ${}_s j^r g(q)\in W$, 
the subspace $E(g,q,W)$ is a $J^r(N)_q$-submodule.
\end{enumerate}
}
\end{definition}

Now, suppose that $P=\mathbb{R}^\ell$. 
The main theorem in \cite{GP} is the following. 
\begin{theorem}[\cite{GP}] \label{@m}
Let $N$ be a manifold of dimension $n$. 
Let $f$ be an embedding 
of $N$ into $\mathbb{R}^m$. 
If $W$ is a modular submanifold of ${}_{s}J^{r}(N,\mathbb{R}^\ell)$ and $m>\ell$, 
then there exists a subset $\Sigma$ with Lebesgue measure zero 
of $\mathcal{L}(\mathbb{R}^{m},\mathbb{R}^{\ell})$ 
such that for any $\pi \in \mathcal{L}(\mathbb{R}^{m},\mathbb{R}^{\ell})-\Sigma$, 
$\pi \circ f:
N\to \mathbb{R}^\ell$ is transverse with respect to $W$.
\end{theorem}
Then, the main theorem in \cite{F1} is the following. 
\begin{theorem}[\cite{F1}]\label{@f}
Let $N$ be a manifold of dimension $n$. 
Let $f$ be an embedding 
of $N$ into an open subset $U$ of $\mathbb{R}^m$. 
Let $F:U\to \mathbb{R}^\ell$ be a mapping. 
If $W$ is a modular submanifold of ${}_{s}J^{r}(N,\mathbb{R}^\ell)$, 
then there exists a subset $\Sigma$ with Lebesgue measure zero 
of $\mathcal{L}(\mathbb{R}^{m},\mathbb{R}^{\ell})$ 
such that for any $\pi \in \mathcal{L}(\mathbb{R}^{m},\mathbb{R}^{\ell})-\Sigma$, 
$F_\pi \circ f:
N\to \mathbb{R}^\ell$ is transverse with respect to $W$.
\end{theorem}
The assertion (6) in Section \ref{section 1}, 
Corollary \ref{embedding} in Section \ref{section 5} and Corollary \ref{stable} in Section \ref{section 6} of the present paper are obtained as corollaries of 
Theorems \ref{main} and \ref{main2} in this paper. 
On the other hand, they are also corollaries of 
Theorem \ref{@f}. 
\section*{Acknowledgements}
The author is most grateful to the anonymous reviewers for carefully reading 
the first manuscript of this paper and for giving invaluable suggestions. 
He is also grateful to Takashi Nishimura for his kind advice, and 
to Atsufumi Honda, Satoshi Koike, 
Osamu 
Saeki, 
Hajime Sato, Masatomo Takahashi, Takahiro Yamamoto for their valuable comments.
The author is supported by JSPS KAKENHI Grant Number 16J06911. 
\>


\begin{thebibliography}{99}
\bibitem{abra}R.~Abraham,  
\textit{Transversality in manifolds of mappings}, 
Bull. Amer. Math. Soc. 
\textbf{69} (1963), 470--474.   



\bibitem{CS}J.~W.~Bruce and P.~J.~Giblin, 
\textit{Curves and singularities $($second edition$)$}, 
Cambridge University Press, Cambridge, 1992. 


\bibitem{GG}M.~Golubitsky and V.~Guillemin, 
\textit{Stable mappings and their singularities}, 
Graduate Texts in Mathematics \textbf{14}, Springer, New York, 1973.



\bibitem{F1}S.~Ichiki, 
\textit{Generic linear perturbations}, 
preprint (available
from arXiv:1607.03220 [math.MG]).


\bibitem{D}S.~Ichiki and T.~Nishimura, 
\textit{Distance-squared mappings}, 
Topology Appl., 
\textbf{160} (2013), 1005--1016.   

\bibitem{L}S.~Ichiki and T.~Nishimura, 
\textit{Recognizable classification of Lorentzian distance-squared 
mappings}, J.~Geom.~Phys., 
\textbf{81} (2014), 62--71.   

\bibitem{G2}S.~Ichiki and T.~Nishimura, 
\textit{Generalized distance-squared mappings of 
$\mathbb{R}^{n+1}$ into $\mathbb{R}^{2n+1}$}, 
Contemporary Mathematics, Amer. Math. Soc., Providence RI, \textbf{675}  (2016),
 121--132.

\bibitem{G3}S.~Ichiki and T.~Nishimura, 
\textit{Preservation of immersed or injective properties by composing generic 
generalized distance-squared mappings}, accepted for publication in Proceedings of Brazil-Mexico 2nd Meeting on Singularities at Salvador, 
to be published in a volume of Springer Proceedings in Mathematics and Statistics.

\bibitem{G1}S.~Ichiki, T.~Nishimura, R.~Oset Sinha and M.~A.~S.~Ruas, 
\textit{Generalized distance-squared mappings of the plane into the plane}, 
Adv. Geom., \textbf{16} (2016), 189--198. 







\bibitem{GP}J. N. Mather, \textit{Generic projections}, 
Ann. of Math., (2), \textbf{98} (1973), 226--245.




%
\end{thebibliography}
\end{document}